\documentclass[11pt]{amsart}
\usepackage[english]{babel}
\usepackage[latin1]{inputenc}

\usepackage[mathcal]{eucal}
\usepackage{amsmath}
\usepackage{amsfonts}
\usepackage{amssymb}
\usepackage{amsthm}
\usepackage{graphicx,epsfig}
\usepackage{amscd}

\DeclareMathOperator{\tr}{tr}

\DeclareMathOperator{\diag}{diag}

\newcommand{\M}{\mathbb{M}}
\newcommand{\R}{\mathbb{R}}
\newcommand{\h}{\mathbb{H}}
\newcommand{\Z}{\mathbb{Z}}

\newcommand{\s}{\mathbb{S}}
\newcommand{\LL}{\mathbb{L}}
\newcommand{\E}{\mathbb{E}}

\newcommand{\psl}{\widetilde{\mathrm{PSL}_2(\mathbb{R})}}

\newcommand{\rmJ}{\mathrm{J}}

\newcommand{\rmT}{\mathrm{T}}
\newcommand{\rmU}{\mathrm{U}}
\newcommand{\rmS}{\mathrm{S}}
\newcommand{\rmd}{\mathrm{d}}

\newcommand{\rmR}{\mathrm{R}}
\newcommand{\rmI}{\mathrm{I}}

\newcommand{\nil}{\mathrm{Nil}_3}

\newcommand{\cM}{{\mathcal M}}

\newcommand{\cF}{{\mathcal F}}
\newcommand{\cE}{{\mathcal E}}
\newcommand{\cV}{{\mathcal V}}

\newcommand{\cH}{{\mathcal H}}
\newcommand{\cA}{{\mathcal A}}

\newcommand{\cB}{{\mathcal B}}

\newcommand{\cD}{{\mathcal D}}

\newcommand{\cZ}{{\mathcal Z}}

\newcommand{\cP}{{\mathcal P}}

\begin{document}

\newtheorem{thm}{Theorem}[section]
\newtheorem*{thmintro}{Theorem}
\newtheorem{cor}[thm]{Corollary}
\newtheorem{prop}[thm]{Proposition}
\newtheorem{app}[thm]{Application}
\newtheorem{lemma}[thm]{Lemma}
\newtheorem{hypothesis}[thm]{Hypothesis}

\theoremstyle{definition}
 
\newtheorem{defn}[thm]{Definition}
\newtheorem*{notation}{Conventions and notations}
\newtheorem{example}[thm]{Example}
\newtheorem{conj}[thm]{Conjecture}
\newtheorem{prob}[thm]{Problem}
 
\theoremstyle{remark}
 
\newtheorem{rem}[thm]{Remark}

\title[Isometric immersions]
{Isometric immersions into \\ $3$-dimensional homogeneous manifolds}
\author{Beno\^\i t Daniel}


\address{IMPA, Estrada Dona Castorina 110,
22460-320 Rio de Janeiro - RJ, BRAZIL}
\email{bdaniel@impa.br}
\urladdr{http://www.math.jussieu.fr/\~{}daniel}

\begin{abstract}
We give a necessary and sufficient condition for a $2$-dimen-sional
Riemannian manifold to be locally 
isometrically immersed into a $3$-dimensional homogeneous manifold with a 
$4$-dimensional isometry group. The condition is expressed
in terms of the metric, the second fundamental form, and data
arising from an ambient Killing field.
This class of $3$-manifolds includes in particular the Berger spheres,
the Heisenberg space $\nil$, the universal cover of the
Lie group $\mathrm{PSL}_2(\R)$ and the product spaces $\s^2\times\R$
and $\h^2\times\R$.
We give some applications to constant mean curvature (CMC)
surfaces in these
manifolds; in particular we prove the existence of a generalized
Lawson correspondence, i.e., a local isometric correspondence between
CMC surfaces in homogeneous $3$-manifolds.
\end{abstract}

\maketitle

\section{Introduction}

A classical problem in geometry is to determine whether a Riemannian
manifold $\cV$ can be isometrically immersed in another Riemaniann
manifold $\bar\cV$. We will restrict ourselves to the case of codimension $1$
immersions, i.e., $\cV$ has dimension $n$ and $\bar\cV$ has 
dimension $n+1$.

It is well known that the Gauss and Codazzi equations 
are necessary
conditions relating the Riemann curvature tensor $\bar\rmR$ of $\bar\cV$,
the Riemann curvature tensor $\rmR$ of $\cV$ and the shape operator
$\rmS$ of $\cV$. Denoting by $\nabla$ the Riemannian connection of 
$\cV$, these equations are the following: 
\begin{equation*}
\langle\rmR(X,Y)Z,W\rangle-\langle\bar\rmR(X,Y)Z,W\rangle
=\langle\rmS X,Z\rangle\langle\rmS Y,W\rangle
-\langle\rmS Y,Z\rangle\langle\rmS X,W\rangle
\end{equation*}
\begin{equation*}
\nabla_X\rmS Y-\nabla_Y\rmS X-\rmS[X,Y]=\bar\rmR(X,Y)N,
\end{equation*}
for all vector fields $X$, $Y$, $Z$ and $W$ on $\cV$.

Moreover, in the case where $\bar\cV$ is a space-form,
i.e., the sphere $\s^{n+1}$, the Euclidean space $\R^{n+1}$ or the
hyperbolic space $\h^{n+1}$, the Gauss and Codazzi equations are also a
sufficient condition for $\cV$ to be locally isometrically immersed
in $\bar\cV$ with $\rmS$ as shape operator. In this case the Gauss and
Codazzi equations involve only the metric and the shape operator of $\cV$.

The author studied this problem when $\bar\cV$ is a product manifold
$\s^n\times\R$ or $\h^n\times\R$ (\cite{codazzi}). Then the Gauss and
Codazzi equations involve the metric of $\cV$, its shape operator $\rmS$, 
the projection $T$ of the vertical vector field (i.e., the unit vector
field corresponding to the factor $\R$) on the tangent space of
$\cV$ and the normal component $\nu$ of the vertical vector field (i.e.,
its inner product with the unit normal of $\cV$). The author proved that
the Gauss and Codazzi equations, together with two other compatibility
equations coming from the fact that the vertical vector field is parallel,
are a necessary and sufficient condition for $\cV$ to be locally
isometrically immersed in $\bar\cV$ with $\rmS$ as shape operator, $T$ as
tangent projection of the vertical vector field and $\nu$ as normal
component of the vertical vector field.

It is natural to try to generalize this result to other homogeneous manifolds.
We will investigate the case of surfaces in manifolds of dimension $3$, i.e.,
$n=2$. Indeed, the classification of simply connected
$3$-dimensional homogeneous manifolds is well known. Such a
manifold has an isometry group of dimension $3$, $4$ or $6$. When the dimension
of the isometry group is $6$, then we have a space form. When the dimension
of the isometry group is $3$, the manifold has the
geometry of the Lie group $\mathrm{Sol}_3$.

In this paper we will consider the homogeneous 
manifolds whose isometry groups have dimension
$4$: such a manifold is a Riemannian fibration over a $2$-dimensional space
form, the fibers are geodesics and there exists a one-parameter family of
translations along the fibers, generated by a unit Killing field $\xi$
which will be called the vertical vector field.
These manifolds are classified, up to isometry, by the curvature $\kappa$ of
the base surface of the fibration
and the bundle curvature $\tau$, where $\kappa$ and $\tau$ can be any real 
numbers satisfying $\kappa\neq 4\tau^2$. The bundle curvature is the number 
$\tau$ such $\bar\nabla_X\xi=\tau X\times\xi$ for any vector field $X$
on $\bar\cV$, where $\bar\nabla$ denotes the Riemannian connection of
$\bar\cV$.

When the bundle curvature $\tau$ vanishes (and then $\kappa\neq 0$),
we get a product manifold
$\M^2(\kappa)\times\R$ where $\M^2(\kappa)$ is the simply connected
$2$-manifold of constant
curvature $\kappa$. Their isometry group has $4$ connected components.
The vertical vector $\xi$ is simply the vector
corresponding to the factor $\R$. This case was treated in
\cite{codazzi}.

When $\tau\neq 0$, the isometry group has $2$ connected components: an 
isometry either preserves the orientations of both the fibers and the 
base of the fibration, or reverses both orientations. These
manifolds are of three types: they have the isometry
group of the Berger spheres for $\kappa>0$, of the Heisenberg space
$\nil$ for $\kappa=0$, and of $\psl$ for
$\kappa<0$. In this paper we will deal with these three types of manifold.
Like for $\M^2(\kappa)\times\R$, the Gauss and Codazzi equations
involve the metric of
$\cV$, its shape operator $\rmS$, the tangential projection $T$ of 
$\xi$ and the normal component $\nu$ of $\xi$. Denoting by $K$ the
curvature of $\rmd s^2$, these equations become
$$K=\det\rmS+\tau^2+(\kappa-4\tau^2)\nu^2,$$
$$\nabla_X\rmS Y-\nabla_Y\rmS X-\rmS[X,Y]=
(\kappa-4\tau^2)\nu(\langle Y,T\rangle X-\langle X,T\rangle Y)$$

The first theorem is the following one.

\begin{thmintro}[theorem \ref{isometry}]
Let $\cV$ be a simply connected oriented
Riemannian manifold of dimension $2$,
$\rmd s^2$ its metric (which we also denote by
$\langle\cdot,\cdot\rangle$), $\nabla$ its  
Riemannian connection and $\rmJ$ the rotation of angle $\frac\pi2$
on $\rmT\cV$. Let $\rmS$ be a field of symmetric operators
$\rmS_y:\rmT_y\cV\to\rmT_y\cV$, $T$ a vector field on $\cV$
and $\nu$ a smooth function on $\cV$ such that
$||T||^2+\nu^2=1$. 

Let $\E$ be a $3$-dimensional homogeneous manifold with a $4$-dimensional
isometry group and $\xi$ its vertical vector field.
Let $\kappa$ be its base curvature and $\tau$ its bundle curvature. 
Then there exists an isometric immersion $f:\cV\to\E$ such that 
the shape operator with respect to the normal $N$ associated to $f$ is
$$\rmd f\circ\rmS\circ\rmd f^{-1}$$ and such that
$$\xi=\rmd f(T)+\nu N$$ 
if and only if $(\rmd s^2,\rmS,T,\nu)$ satisfies the
Gauss and Codazzi equations for $\E$ and, for all vector fields
$X$ on $\cV$, the following equations:
$$\nabla_XT=\nu(\rmS X-\tau\rmJ X),\quad
\rmd\nu(X)+\langle\rmS X-\tau\rmJ X,T\rangle=0.$$
In this case, the
immersion is unique up to a global isometry of $\E$ preserving
the orientations of both the fibers and the base of the fibration.
\end{thmintro}

The two additional conditions come from the fact that
$\bar\nabla_X\xi=\tau X\times\xi$ for all vector fields $X$.

We notice that this theorem seems specific to dimension $2$, since the
operator of rotation $\rmJ$ is involved.

The method to prove this theorem is similar to that of \cite{codazzi}
and was inspired by that of Tenenblat
(\cite{tenenblat}): it is based on differential forms, moving frames
and integrable distributions. 
However, things are technically much more complicated
here: in \cite{codazzi} the proof was simplified by the fact that
$\s^n\times\R$ and $\h^n\times\R$ can be included in $\R^{n+2}$ and
in the Lorentz space $\LL^{n+2}$ respectively. We will first present
the models used for the $3$-dimensional homogeneous manifolds, and then we will prove
the theorem.

Finally we will give two applications of the main theorem
to constant mean curvature (CMC)
surfaces in $3$-dimensional homogeneous manifolds
with $4$-dimensional isometry group.

The first application is the existence of an isometric correspondence 
between certain CMC surfaces in homogeneous $3$-manifolds with the same
anisotropy coefficient $\kappa-4\tau^2$. This correspondence
generalizes the classical Lawson correspondence between certain
CMC surfaces in space-forms. This is the following theorem.

\begin{thmintro}[see theorem \ref{sisters}]
Let $\E_1$ and $\E_2$ be two $3$-dimensional homogeneous manifolds
with $4$-dimensional isometry groups, of base curvatures $\kappa_1$
and $\kappa_2$ and bundle curvatures $\tau_1$ and $\tau_2$
respectively, and such that
$$\kappa_1-4\tau_1^2=\kappa_2-4\tau_2^2.$$

Let $H_1$ and $H_2$ be two real numbers such that
$$\tau_1^2+H_1^2=\tau_2^2+H_2^2.$$ 

Then there exists an isometric 
correspondence between simply connected 
CMC $H_1$ surfaces in $\E_1$ and
simply connected
CMC $H_2$ surfaces in $\E_2$.

This correspondence is called the correspondence of the sister
surfaces.
\end{thmintro}

The second application is the existence of ``twin immersions'' of
non-minimal CMC immersions in homogeneous $3$-manifolds with non-vanishing
bundle curvature. This twin immersion might be useful to prove an
Alexan-drov-type theorem in these manifolds.

\begin{notation}
In this paper we will use the following index conventions: Latin
letters $i$, $j$, etc, denote integers between $1$ and $n$ (or
the integers $1$ and $2$), Greek
letters $\alpha$, $\beta$, etc, denote integers between $1$ and
$n+1$ (or between $1$ and $3$). 

The set of vector fields on a Riemannian manifold $\cV$ will be
denoted by $\mathfrak{X}(\cV)$.

The Riemann curvature tensor $\rmR$ of a Riemannian manifold $\cV$ 
of Riemannian connection $\nabla$ is
defined using the following convention:
$$\rmR(X,Y)Z=\nabla_Y\nabla_XZ-\nabla_X\nabla_YZ+\nabla_{[X,Y]}Z.$$

The shape operator of a hypersurface $\cV$ of a Riemannian manifold
$\bar\cV$ associated to its unit normal $N$ is
$$\rmS X=-\bar\nabla_XN$$ where $\bar\nabla$ is the Riemannian
connection of $\bar\cV$. 
\end{notation}

\section{$3$-dimensional homogeneous manifolds with $4$-dimensional isometry group}

In this section we will give the general setting for 
simply connected homogeneous 
$3$-manifolds with $4$-dimensional isometry group and we will describe the 
models used. We will consider only those having non-vanishing bundle curvature
(since the product manifolds $\M^2(\kappa)\times\R$ were treated in
\cite{codazzi}).
The reader can refer to \cite{scott} for the geometry of
$3$-dimensional homogeneous manifolds.

\subsection{Canonical frame} \label{canonicalframe}

Let $\E$ be a simply connected
$3$-dimensional homogeneous manifold with a $4$-dimensional isometry group.
Such a manifold is a Riemannian fibration over a 
simply connected $2$-manifold of constant
curvature $\kappa$. The fibers are geodesics. We will denote by
$\xi$ a unit vector field on $\E$ tangent to the fibers; it will be called
the vertical vector field. It is a Killing field (corresponding to
translations along the fibers).

We will denote by $\bar\nabla$ and $\bar\rmR$ the Riemannian connection
and the Riemannian curvature tensor of $\E$ respectively.

We assume that $\E$ is not a product manifold $\M^2(\kappa)\times\R$.

The manifold $\E$ locally has a direct orthonormal
frame $(E_1,E_2,E_3)$ with $$E_3=\xi$$
whose non-vanishing Christoffel symbols 
$\bar\Gamma^\alpha_{\beta\gamma}
=\langle\nabla_{E_\alpha}E_\beta,E_\gamma\rangle$ are the following:
$$\bar\Gamma^3_{12}=\bar\Gamma^1_{23}=-\bar\Gamma^3_{21}
=-\bar\Gamma^2_{13}=\tau,$$
$$\bar\Gamma^1_{32}=-\bar\Gamma^2_{31}=\tau-\sigma,$$
for some real numbers $\sigma$ and $\tau\neq 0$ 
(this will be explicited in the sequel).
Then we have
$$[E_1,E_2]=2\tau E_3,\quad
[E_2,E_3]=\sigma E_1,\quad
[E_3,E_1]=\sigma E_2.$$
We will call $(E_1,E_2,E_3)$ the canonical frame of $\E$.
For all vector field $X$ we have $$\bar\nabla_XE_3=\tau X\times E_3$$ 
where $\times$ denotes the vector product in $\E$, i.e., for all
vector fields $X$, $Y$, $Z$,
$\langle X\times Y,Z\rangle=\det_{(E_1,E_2,E_3)}(X,Y,Z)$.

Setting $$\langle\bar\rmR(X\wedge Y),Z\wedge W\rangle
=\langle\bar\rmR(X,Y)Z,W\rangle,$$ 
the matrix of $\bar\rmR$ in the basis 
$(E_2\wedge E_3,E_3\wedge E_1,E_1\wedge E_2)$ is
$$\bar\rmR=\diag(a,a,b)$$ with $$a=\tau^2,\quad
b=-3\tau^2+2\sigma\tau.$$ 

%
%

We now compute the curvature $\kappa$ 
of the base of the fibration. If $\bar M\to M$ is a Riemannian submersion,
then the sectional curvature of a $2$-plane $\Pi$ in $M$ generated by
an orthonormal pair $(X,Y)$ is
$$K(\Pi)=\bar K(\bar\Pi)+
\frac34\left|\left|[\bar X,\bar Y]^{\mathrm v}\right|\right|^2$$
where $\bar X$ and $\bar Y$ are horizontal lifts of $X$ and $Y$ in $\bar M$,
$\bar K(\bar\Pi)$ is the sectional curvature of a $2$-plane $\bar\Pi$
in $\bar M$ generated by $(\bar X,\bar Y)$, and where $Z^{\mathrm v}$
denotes the vertical part of a vector field $Z$ in $\bar M$
(see \cite{docarmo}, chapter 8).
In our case we get
$$\kappa=\langle\bar\rmR(E_1,E_2)E_1,E_2\rangle
+\frac34\left|\left|[E_1,E_2]^{\mathrm v}\right|\right|^2
=b+\frac34\left|\left|2\tau E_3^{\mathrm v}\right|\right|^2
=b+3\tau^2.$$
Thus we have $b=\kappa-3\tau^2$, and so
$$\sigma=\frac\kappa{2\tau}.$$

\begin{prop} \label{exprbarR}
For all vector fields $X,Y,Z,W$ on $\E$ we have
$$\langle\bar\rmR(X,Y)Z,W\rangle=
(\kappa-3\tau^2)\langle \rmR_0(X,Y)Z,W\rangle
+(\kappa-4\tau^2)\langle\rmR_1(\xi;X,Y)Z,W\rangle$$ 
with $$\rmR_0(X,Y)Z=\langle X,Z\rangle Y-\langle Y,Z\rangle X,$$
\begin{eqnarray*}
\rmR_1(V;X,Y)Z & = &
\langle Y,V\rangle\langle Z,V\rangle X
+\langle Y,Z\rangle\langle X,V\rangle V \\
& & -\langle X,Z\rangle\langle Y,V\rangle V
-\langle X,V\rangle\langle Z,V\rangle Y.
\end{eqnarray*}
\end{prop}

\begin{proof}
We set $X=\tilde X+x\xi$ with $\tilde X$ horizontal and $x=\langle
X,\xi\rangle$, etc. Using the multilinearity
of the Riemann curvature tensor, we get a sum of 16 terms; the terms
where $\xi$ appears three or four times, or twice at positions $1,2$
or $3,4$, vanish by antisymmetry. The terms where $\xi$ appears once
vanish because the matrix of $\bar\rmR$ in the basis
$(E_2\wedge E_3,E_3\wedge E_1,E_1\wedge E_2)$ is diagonal. Hence we have
\begin{eqnarray*}
\langle\bar\rmR(X,Y)Z,W\rangle & = &
\langle\bar\rmR(\tilde X,\tilde Y)\tilde Z,\tilde W\rangle \\
& & +yw\langle\bar\rmR(\tilde X,\xi)\tilde Z,\xi\rangle
+yz\langle\bar\rmR(\tilde X,\xi)\xi,\tilde W\rangle \\
& & +xw\langle\bar\rmR(\xi,\tilde Y)\tilde Z,\xi\rangle
+xz\langle\bar\rmR(\xi,\tilde Y)\xi,\tilde W\rangle \\
& = & (\kappa-3\tau^2)(
\langle\tilde X,\tilde Z\rangle\langle\tilde Y,\tilde W\rangle
-\langle\tilde X,\tilde W\rangle\langle\tilde Y,\tilde Z\rangle) \\
& & +\tau^2(yw\langle\tilde X,\tilde Z\rangle
-yz\langle\tilde X,\tilde W\rangle
-xw\langle\tilde Y,\tilde Z\rangle
+xz\langle\tilde Y,\tilde W\rangle) \\
& = & (\kappa-3\tau^2)(
\langle X,Z\rangle\langle Y,W\rangle
-\langle X,W\rangle\langle Y,Z\rangle) \\
& & -(\kappa-4\tau^2)
(\langle X,Z\rangle\langle Y,\xi\rangle\langle W,\xi\rangle
+\langle Y,W\rangle\langle X,\xi\rangle\langle Z,\xi\rangle \\
& & \quad-\langle X,W\rangle\langle Y,\xi\rangle\langle Z,\xi\rangle
-\langle Y,Z\rangle\langle X,\xi\rangle\langle W,\xi\rangle).
\end{eqnarray*}
\end{proof}

%

\subsection{The manifolds with the isometry group of the Berger spheres}
\label{bergerspheres}

They occur when $\tau\neq 0$ and $\kappa>0$; they are
fibrations over round $2$-spheres.
They are obtained by deforming the metric of a round sphere in a way
preserving the Hopf fibration but modifying the length of the fibers.
Their isometry group is included in that of the round sphere. The reader
can refer to \cite{petersen}.

The sphere $\s^3$ is the univeral covering of $\mathrm{SO}_3(\R)$, which
can be identified with the unitary
tangent bundle to the $2$-sphere $\rmU\s^2$. 
Indeed, the group $\mathrm{SO}_3(\R)$ acts
transitively on $\rmU\s^2$, and the stabilizer of any point in $\rmU\s^2$
is trivial. The unitary tangent bundle $\rmU\s^2$ can be endowed 
with the metric induced by
the standard metric on the tangent bundle $\rmT\s^2$. We will give an
expression of this metric.

Let $(x,y)\mapsto\varphi(x,y)$ be a conformal parametrization 
of a domain $D$ in $\s^2$
and let $\lambda$ be the conformal factor, i.e., the metric of $D$ is
given by $\lambda^2(\rmd x^2+\rmd y^2)$. Then a parametrization of 
$\rmU D$ is the following:
$$(x,y,\theta)\mapsto
\left(\varphi(x,y),\frac1\lambda(\cos\theta\partial_x
+\sin\theta\partial_y)\right).$$

Let $p=\varphi(x,y)\in D$, $v\in\rmT_pD$ and
$V\in\rmT_{(p,v)}(\rmU D)$. Let $\alpha(t)=(p(t),v(t))$ be
a curve such that $v(t)\in\rmT_{p(t)}\h^2$, $p(0)=p$, $v(0)=v$ and
$\alpha'(0)=V$. Then the norm of $V$ is given by
$$||V||_{(p,v)}^2=||\rmd\pi(V)||_p^2+
\left|\left|\frac{\mathrm{D}v}{\rmd t}(0)\right|\right|_p^2$$
where $\pi:\rmU D\to D$ is the canonical projection.

We set $\alpha(t)=(x(t),y(t),\theta(t))$. Then we have
$$v(t)=\frac1\lambda(\cos\theta(t)\partial_x+\sin\theta(t)\partial_y),$$
and thus
\begin{eqnarray*}
\frac{\mathrm{D}v}{\rmd t} & = & 
-\frac{\dot\lambda}{\lambda^2}
(\cos\theta\partial_x+\sin\theta\partial_y)
+\frac{\dot\theta}\lambda(-\sin\theta\partial_x
+\cos\theta\partial_y) \\
& & +\frac1\lambda(\cos\theta(\dot x\nabla_{\partial_x}\partial_x
+\dot y\nabla_{\partial_y}\partial_x)
+\sin\theta(\dot x\nabla_{\partial_x}\partial_y
+\dot y\nabla_{\partial_y}\partial_y)),
\end{eqnarray*}
where the dot denotes the derivation with respect to $t$.
Since $\dot\lambda=\dot x\lambda_x+\dot y\lambda_y$,
$\nabla_{\partial_x}\partial_x=\frac{\lambda_x}\lambda\partial_x
-\frac{\lambda_y}\lambda\partial_y$,
$\nabla_{\partial_y}\partial_y=-\frac{\lambda_x}\lambda\partial_x
+\frac{\lambda_y}\lambda\partial_y$
and $\nabla_{\partial_x}\partial_y=
\nabla_{\partial_y}\partial_x=\frac{\lambda_y}\lambda\partial_x
+\frac{\lambda_x}\lambda\partial_y$, we get
$$\frac{\mathrm{D}v}{\rmd t}=\frac1{\lambda^2}
(\lambda\dot\theta+\dot y\lambda_x-\dot x\lambda_y)
(\cos\theta\partial_y-\sin\theta\partial_x).$$

Thus $$||V||^2_{(p,v)}=\lambda^2(\dot x^2+\dot y^2)
+\frac1{\lambda^2}(\lambda\dot\theta+\dot y\lambda_x-\dot x\lambda_y)^2.$$

Setting $z=\theta$ on the universal covering, we get the following
expression for the metric of $\widetilde{\rmU D}$:
$$\rmd s^2=\lambda^2(\rmd x^2+\rmd y^2)
+\left(-\frac{\lambda_y}\lambda\rmd x
+\frac{\lambda_x}\lambda\rmd y+\rmd z\right)^2.$$

We now choose $D=\s^2\setminus\{\infty\}$ with the metric of constant
curvature $4$ (i.e., the metric of the round sphere of radius $\frac12$)
given by the stereographic projection, i.e.,
$$\lambda=\frac 1{1+x^2+y^2}.$$ Then we get
$$\rmd s^2=\lambda^2(\rmd x^2+\rmd y^2)
+(2\lambda(y\rmd x-x\rmd y)+\rmd z)^2.$$

More generally, $\R^3$ endowed with the metric
$$\rmd s^2=
\lambda^2(\rmd x^2+\rmd y^2)
+\left(\tau\lambda(y\rmd x-x\rmd y)+\rmd z\right)^2$$ with
$$\lambda=\frac1{1+\frac\kappa4(x^2+y^2)}$$
is the universal cover of a homogeneous manifold $\E$ of bundle curvature $\tau$ 
and of base curvature $\kappa>0$ minus the fiber
corresponding to the point $\infty\in\s^2$. 
The fibers are given by $\{x=x_0,y=y_0\}$ in these coordinates.
The canonical frame is $(E_1,E_2,E_3)$ with
\begin{equation} \label{canonicalbergerspheres}
\begin{array}{c}
E_1=\lambda^{-1}(\cos(\sigma z)\partial_x+\sin(\sigma z)\partial_y)
+\tau(x\sin(\sigma z)-y\cos(\sigma z))\partial_z, \\
E_2=\lambda^{-1}(-\sin(\sigma z)\partial_x+\cos(\sigma z)\partial_y)
+\tau(x\cos(\sigma z)+y\sin(\sigma z))\partial_z, \\
E_3=\partial_z
\end{array}
\end{equation}
with $$\sigma=\frac\kappa{2\tau},$$
which satisfy
$$[E_1,E_2]=2\tau E_3,\quad [E_2,E_3]=\frac\kappa{2\tau}E_1,
\quad [E_3,E_1]=\frac\kappa{2\tau}E_2.$$
This frame is defined on the open set $\E'$ which is $\E$ minus the fiber
corresponding to the point $\infty\in\s^2$.

%


The Berger spheres in the strict sense are the manifolds such that $\kappa=4$.

\subsection{The manifolds with the isometry group of the
Heisenberg space $\nil$}
\label{heisenberg}

They occur when $\tau\neq 0$ and $\kappa=0$; they are fibrations
over the Euclidean plane.

The Heisenberg space is the Lie group
$$\nil=\left\{\left(\begin{array}{ccc}
1 & a & c \\
0 & 1 & b \\
0 & 0 & 1
\end{array}\right);(a,b,c)\in\R^3\right\}$$
endowed with a left invariant metric.

It is useful to use exponential coordinates. In this model, the Heisenberg
space $\nil$ is $\R^3$ endowed with the following metric:
$$\rmd s^2=\rmd x^2+\rmd y^2+
(\tau(y\rmd x-x\rmd y)+\rmd z)^2.$$
The fibers are given by $\{x=x_0,y=y_0\}$ in these coordinates.

The canonical frame is $(E_1,E_2,E_3)$ with
\begin{equation} \label{canonicalheisenberg}
E_1=\partial_x-\tau y\partial_z,\quad
E_2=\partial_y+\tau x\partial_z,\quad
E_3=\partial_z,
\end{equation}
which satisfy
$$[E_1,E_2]=2\tau E_3,\quad [E_2,E_3]=0,
\quad [E_3,E_1]=0.$$

The reader can refer to \cite{mercuri} (where $\tau=\frac12$).

\subsection{The manifolds with the isometry group of $\psl$}
\label{psl}




They occur when $\tau\neq 0$ and $\kappa<0$; they are fibrations over
hyperbolic planes.

The Lie group $\psl$ with its standard metric
can be identified with the universal covering of the unitary
tangent bundle to the hyperbolic plane $\rmU\h^2$ equipped with its
canonical metric. Indeed, the group $\mathrm{PSL}_2(\R)$ acts
transitively on $\rmU\h^2$, and the stabilizer of any point 
in $\rmU\h^2$ is trivial. The unitary tangent bundle $\rmU\h^2$
can be endowed with the metric induced by
the standard metric on the tangent bundle $\rmT\h^2$. 
The reader can refer to \cite{scott}. We will give an
expression of this metric. 

Let $(x,y)\mapsto\varphi(x,y)$ be a conformal parametrization of $\h^2$
and let $\lambda$ be the conformal factor, i.e., the metric of $\h^2$ is
given by $\lambda^2(\rmd x^2+\rmd y^2)$. Then, proceeding as in section
\ref{bergerspheres}, we obtain that a metric on $\psl$ is
$$\rmd s^2=\lambda^2(\rmd x^2+\rmd y^2)
+\left(-\frac{\lambda_y}\lambda\rmd x
+\frac{\lambda_x}\lambda\rmd y+\rmd z\right)^2.$$
This metric defines a homogeneous manifold with $\kappa=-1$ and 
$\tau=-\frac12$.

%

%
%

More generally, we can take the Poincar\'e disk model for the hyperbolic
 plane of constant curvature $\kappa<0$. The manifold 
$\mathbb{D}^2\left(\frac2{\sqrt{-\kappa}}\right)\times\R$, where
$\mathbb{D}^2(\rho)=\{(x,y)\in\R^2;x^2+y^2<\rho^2\}$,
endowed with the metric
$$\rmd s^2=
\lambda^2(\rmd x^2+\rmd y^2)
+\left(\tau\lambda(y\rmd x-x\rmd y)+\rmd z\right)^2$$ with
$$\lambda=\frac1{1+\frac\kappa4(x^2+y^2)}$$
is a homogeneous manifold of bundle curvature $\tau$ and of base curvature
$\kappa<0$. 
The fibers are given by $\{x=x_0,y=y_0\}$ in these coordinates.
The canonical frame
is $(E_1,E_2,E_3)$ with
\begin{equation} \label{canonicalpsl}
\begin{array}{c}
E_1=\lambda^{-1}(\cos(\sigma z)\partial_x+\sin(\sigma z)\partial_y)
+\tau(x\sin(\sigma z)-y\cos(\sigma z))\partial_z, \\
E_2=\lambda^{-1}(-\sin(\sigma z)\partial_x+\cos(\sigma z)\partial_y)
+\tau(x\cos(\sigma z)+y\sin(\sigma z))\partial_z, \\
E_3=\partial_z
\end{array}
\end{equation}
with $$\sigma=\frac\kappa{2\tau},$$
which satisfy
$$[E_1,E_2]=2\tau E_3,\quad [E_2,E_3]=\frac\kappa{2\tau}E_1,
\quad [E_3,E_1]=\frac\kappa{2\tau}E_2.$$

\section{Preliminaries}

\subsection{The compatibility equations for surfaces in
$3$-dimensional homogeneous manifolds} \label{compatibilityE}

We consider a $3$-dimensional homogeneous manifold $\E$ with an
isometry group of dimension $4$, of bundle curvature $\tau$ and of base
curvature $\kappa$. Let $\bar\rmR$ be the Riemann curvature tensor of $\E$.
Let $\cV$ be an oriented surface in $\E$, 
$\nabla$ the Riemannian connection of $\cV$,
$\rmJ$ the rotation of angle $\frac\pi2$ on $\rmT\cV$,
$N$ the unit normal to $\cV$
and $\rmS$ the shape operator of $\cV$.

\begin{prop}
For $X,Y,Z,W\in\mathfrak{X}(\cV)$ we have
$$\langle\bar\rmR(X,Y)Z,W\rangle=
(\kappa-3\tau^2)\langle \rmR_0(X,Y)Z,W\rangle
+(\kappa-4\tau^2)\langle\rmR_1(T;X,Y)Z,W\rangle,$$ 
$$\bar\rmR(X,Y)N=(\kappa-4\tau^2)\nu
(\langle Y,T\rangle X-\langle X,T\rangle Y),$$
where $$\nu=\langle N,\xi\rangle,$$
$T$ is the projection of $\xi$ on 
$\rmT\cV$, i.e., $$T=\xi-\nu N,$$
and $\rmR_0$ and $\rmR_1$ are as in proposition \ref{exprbarR}.
\end{prop}

\begin{proof}
This is a consequence of proposition \ref{exprbarR}, using the fact
that $X$, $Y$ and $Z$ are tangent to the
surface and $N$ is normal to the surface.
\end{proof}

\begin{cor}
The Gauss and Codazzi equations in $\E$ are
$$K=\det\rmS+\tau^2+(\kappa-4\tau^2)\nu^2,$$
$$\nabla_X\rmS Y-\nabla_Y\rmS X-\rmS[X,Y]=
(\kappa-4\tau^2)\nu(\langle Y,T\rangle X-\langle X,T\rangle Y),$$
where $K$ is the Gauss curvature of $\cV$.
\end{cor}

\begin{prop}
For $X\in\mathfrak{X}(\cV)$ we have
$$\nabla_XT=\nu(\rmS X-\tau\rmJ X),\quad
\rmd\nu(X)+\langle\rmS X-\tau\rmJ X,T\rangle=0.$$
\end{prop}

\begin{proof}
On the one hand we have
\begin{eqnarray*}
\bar\nabla_X\xi & = & \bar\nabla_X(T+\nu N) \\
& = & \bar\nabla_XT+\rmd\nu(X)N+\nu\bar\nabla_XN \\
& = & \nabla_XT+\langle\rmS X,T\rangle N+\rmd\nu(X)N-\nu\rmS X.
\end{eqnarray*}

On the other hand we have
\begin{eqnarray*}
\bar\nabla_X\xi & = & \tau X\times\xi \\
& = & \tau X\times(T+\nu N) \\
& = & \tau(\langle\rmJ X,T\rangle N-\nu\rmJ X).
\end{eqnarray*}

We conclude taking the tangential and normal parts in both expressions.
\end{proof}

\subsection{Moving frames} \label{movingframes}

In this section we introduce some material about the technique of
moving frames.

Let $\cV$ be a Riemannian manifold of dimension $n$, 
$\nabla$ its Levi-Civita connection, and
$\rmR$ the Riemannian curvature tensor. 
Let $\rmS$ be a field of symmetric operators
$\rmS_y:\rmT_y\cV\to\rmT_y\cV$.
Let $(e_1,\dots,e_n)$ be a local orthonormal frame on $\cV$ and
$(\omega^1,\dots,\omega^n)$ the dual basis of  $(e_1,\dots,e_n)$,
i.e., $$\omega^i(e_k)=\delta^i_k.$$ We also set
$$\omega^{n+1}=0.$$ 

We define the forms $\omega^i_j$, $\omega^{n+1}_j$, $\omega^i_{n+1}$
and $\omega^{n+1}_{n+1}$ on $\cV$ by
$$\omega^i_j(e_k)=\langle\nabla_{e_k}e_j,e_i\rangle,\quad
\omega^{n+1}_j(e_k)=\langle\rmS e_k,e_j\rangle,$$
$$\omega^j_{n+1}=-\omega^{n+1}_j,\quad
\omega^{n+1}_{n+1}=0.$$
Then we have
$$\nabla_{e_k}e_j=\sum_i\omega^i_j(e_k)e_i,\quad
\rmS e_k=\sum_j\omega^{n+1}_j(e_k)e_j.$$

Finally we set $R^i_{klj}=\langle\rmR(e_k,e_l)e_j,e_i\rangle$.

\begin{prop} \label{differentiation}
We have the following formulas:
\begin{equation} \label{diffomega1}
\rmd\omega^i+\sum_p\omega^i_p\wedge\omega^p=0,
\end{equation}
\begin{equation} \label{diffomega2}
\sum_p\omega^{n+1}_p\wedge\omega^p=0,
\end{equation}
\begin{equation} \label{diffomega3}
\rmd\omega^i_j+\sum_p\omega^i_p\wedge\omega^p_j=
-\frac{1}{2}\sum_k\sum_lR^i_{klj}\omega^k\wedge\omega^l,
\end{equation}
\begin{equation} \label{diffomega4}
\rmd\omega^{n+1}_j+\sum_p\omega^{n+1}_p\wedge\omega^p_j=
\frac{1}{2}\sum_k\sum_l\langle\nabla_{e_k}\rmS e_l
-\nabla_{e_l}\rmS e_k-\rmS[e_k,e_l],e_j\rangle\omega^k\wedge\omega^l.
\end{equation}
\end{prop}

For a proof of these classical formulas, 
the reader can refer to \cite{codazzi}, proposition 2.4.

\subsection{Some facts about hypersurfaces}
\label{hypersurfaces}

In this section we consider an orientable hypersurface $\cV$ of an
$(n+1)$-dimensionnal Riemannian manifold $\bar\cV$.
Let $(e_1,\dots,e_n)$ be a local 
orthonormal frame on $\cV$, $e_{n+1}$ the normal to $\cV$, and
$(E_1,\dots,E_{n+1})$ a local orthonormal frame on $\bar\cV$. We denote
by $\nabla$ and $\bar\nabla$ the Riemannian connections on
$\cV$ and $\bar\cV$ respectively, and by $\rmS$ the shape operator of
$\cV$ (with respect to the normal $e_{n+1}$). We define the forms
$\omega^\alpha$, $\omega^\alpha_\beta$ on $\cV$ as in
section \ref{movingframes}. Then we have
$$\bar\nabla_{e_k}e_\beta=\sum_\gamma\omega^\gamma_\beta(e_k)e_\gamma.$$

Let $A\in\mathrm{SO}_{n+1}(\R)$ be the matrix whose columns are the
coordinates of the $e_\beta$ in the frame $(E_\alpha)$, namely
$A^\alpha_\beta=\langle e_\beta,E_\alpha\rangle$. Let
$\Omega=(\omega^\alpha_\beta)\in\cM_{n+1}(\R)$. 

\begin{lemma} \label{diffA}
The matrix $A$ satisfies the following equation:
$$A^{-1}\rmd A=\Omega+L(A)$$
with $$L(A)^\alpha_\beta=\sum_k
\left(\sum_{\gamma,\delta,\varepsilon}A^\varepsilon_\alpha
A^\gamma_kA^\delta_\beta
\bar\Gamma_{\gamma\varepsilon}^\delta\right)\omega^k,$$
where the $\bar\Gamma_{\gamma\varepsilon}^\delta$ are the Christoffel
symbols of the frame $(E_\alpha)$.
\end{lemma}

\begin{proof}
We have $$e_\beta=\sum_\alpha A^\alpha_\beta E_\alpha.$$
Then, on the one hand we have
\begin{eqnarray*}
\bar\nabla_{e_k}e_\beta
& = & \sum_\delta\rmd A^\delta_\beta(e_k)E_\delta
+\sum_\delta A^\delta_\beta\bar\nabla_{e_k}E_\delta \\
& = & \sum_\varepsilon\rmd A^\varepsilon_\beta(e_k)E_\delta
+\sum_\gamma\sum_\delta\sum_\varepsilon
A^\delta_\beta A^\gamma_k\bar\Gamma^\varepsilon_{\gamma\delta}
E_\varepsilon,
\end{eqnarray*}
and on the other hand we have
$$\bar\nabla_{e_k}e_\beta=
\sum_\gamma\sum_\varepsilon\omega^\gamma_\beta(e_k)
A^\varepsilon_\gamma E_\varepsilon.$$ 
Identifying the coefficients we get
\begin{eqnarray*}
\rmd A^\varepsilon_\beta(e_k) & = &
-\sum_\gamma\sum_\delta
A^\delta_\beta A^\gamma_k\bar\Gamma^\varepsilon_{\gamma\delta}
+\sum_\gamma\omega^\gamma_\beta(e_k)A^\varepsilon_\gamma \\
& = & \sum_\gamma\sum_\delta
A^\delta_\beta A^\gamma_k\bar\Gamma^\delta_{\gamma\varepsilon}
+\sum_\gamma\omega^\gamma_\beta(e_k)A^\varepsilon_\gamma
\end{eqnarray*}
since the frame $(E_\alpha)$ is orthonormal.

We conclude using the fact that $A^{-1}$ is the transpose of $A$.
\end{proof}

\section{Isometric immersions of surfaces into $3$-dimensional 
homogeneous manifolds}

We consider a simply connected oriented
Riemannian manifold $\cV$ of dimension
$2$. Let $\rmd s^2$ be the metric on $\cV$ (we will also denote it by
$\langle\cdot,\cdot\rangle$), $\nabla$  
the Riemannian connection of $\cV$, $\rmR$ its Riemann curvature
tensor and $\rmJ$ the rotation of angle $\frac\pi2$ on $\rmT\cV$.
Let $\rmS$ be a field of symmetric operators
$\rmS_y:\rmT_y\cV\to\rmT_y\cV$, $T$ a vector field on $\cV$
such that $||T||\leqslant 1$ and $\nu$ a smooth function on $\cV$ 
such that $\nu^2\leqslant 1$.

The compatibility equations for surfaces in $3$-dimensional
homogeneous manifolds with $4$-dimensional isometry group
established in section \ref{compatibilityE} suggest to introduce
the following definition.

\begin{defn}
Let $\E$ be a $3$-dimensional homogeneous manifold with a $4$-dimensional
isometry group. Let $\kappa$ be its base curvature and $\tau$ its 
bundle curvature.
We say that $(\rmd s^2,\rmS,T,\nu)$ satisfies the
compatibility equations for $\E$
if $$||T||^2+\nu^2=1$$ and, for all $X,Y,Z\in\mathfrak{X}(\cV)$, 
\begin{equation} \label{gaussE}
K=\det\rmS+\tau^2+(\kappa-4\tau^2)\nu^2,
\end{equation}
\begin{equation} \label{codazziE}
\nabla_X\rmS Y-\nabla_Y\rmS X-\rmS[X,Y]=
(\kappa-4\tau^2)\nu(\langle Y,T\rangle X-\langle X,T\rangle Y),
\end{equation}
\begin{equation} \label{conditionT1}
\nabla_XT=\nu(\rmS X-\tau\rmJ X),
\end{equation}
\begin{equation} \label{conditionT2}
\rmd\nu(X)+\langle\rmS X-\tau\rmJ X,T\rangle=0.
\end{equation}
\end{defn}

\begin{rem}
We notice that \eqref{conditionT1} implies \eqref{conditionT2} except
when $\nu=0$ (by differentiating the identity $\langle T,T\rangle
+\nu^2=1$ with respect to $X$).
\end{rem}

\begin{thm} \label{isometry}
Let $\cV$ be a simply connected oriented
Riemannian manifold of dimension $2$,
$\rmd s^2$ its metric and $\nabla$ its Riemannian connection.
Let $\rmS$ be a field of symmetric operators
$\rmS_y:\rmT_y\cV\to\rmT_y\cV$, $T$ a vector field on $\cV$ and $\nu$
a smooth function on $\cV$ such that $||T||^2+\nu^2=1$. 

Let $\E$ be a $3$-dimensional homogeneous manifold with a $4$-dimensional
isometry group
and $\xi$ its vertical vector field.
Let $\kappa$ be its base curvature and $\tau$ its bundle curvature.
Then there exists an isometric immersion $f:\cV\to\E$
such that
the shape operator with respect to the normal $N$ associated to $f$ is
$$\rmd f\circ\rmS\circ\rmd f^{-1}$$ and such that
$$\xi=\rmd f(T)+\nu N$$
if and only if $(\rmd s^2,\rmS,T,\nu)$ satisfies the
compatibility equations for $\E$.
In this case, the
immersion is unique up to a global isometry of $\E$ preserving
the orientations of both the fibers and the base of the fibration.
\end{thm}

The fact that the compatibility equations are necessary was proved in
section \ref{compatibilityE}. 
To prove that they are sufficient, we consider a local orthonormal frame
$(e_1,e_2)$ on $\cV$ and the forms $\omega^i$, $\omega^3$,
$\omega^i_j$, $\omega^3_j$, $\omega^i_3$ and
$\omega^3_3$ as in section \ref{movingframes} (with $n=2$).

From now on we assume that $\tau\neq 0$ since the case $\tau=0$
was treated in \cite{codazzi}.

We denote by 
$(E_1,E_2,E_3)$ the canonical frame of $\E$ (see section 
\ref{canonicalframe}); in particular we
have $E_3=\xi$. We denote by $\E'$ the open set where
the canonical frame is defined (in particular we have $\E'=\E$ when
$\kappa=0$ or $\kappa<0$; see sections \ref{bergerspheres},
\ref{heisenberg} and \ref{psl}).

We set 
$$T^k=\langle T,e_k\rangle,\quad T^3=\nu.$$
We define the one-form $\eta$ on $\cV$ by
$$\eta(X)=\langle T,X\rangle.$$
In the frame $(e_1,e_2)$ we have $\eta=\sum_kT^k\omega^k$.
We define the following matrix of one-forms:
$$\Omega=(\omega^\alpha_\beta)\in\cM_3(\R).$$
For $Z\in\mathrm{SO}_3(\R)$, we set
$$L(Z)^\alpha_\beta=\sum_k
\left(\sum_{\gamma,\delta,\varepsilon}Z^\varepsilon_\alpha
Z^\gamma_kZ^\delta_\beta
\bar\Gamma_{\gamma\varepsilon}^\delta\right)\omega^k,$$
where the $\bar\Gamma_{\gamma\varepsilon}^\delta$ are the Christoffel
symbols of the frame $(E_\alpha)$ (see section \ref{hypersurfaces}).
This defines an antisymmetric matrix of $1$-forms.

We also set $\sigma=\frac\kappa{2\tau}$.

From now on we assume that the hypotheses of theorem \ref{isometry}
are satisfied. We first prove some technical lemmas that are
consequences of the compatibility equations.

\begin{lemma} \label{diffeta}
We have $$\rmd\eta=-2\tau\nu\omega^1\wedge\omega^2.$$
\end{lemma}

\begin{proof}
By \eqref{conditionT1} we have
$\rmd\eta(X,Y)=\langle\nabla_XT,Y\rangle-\langle\nabla_YT,X\rangle
=2\tau\nu\langle X,\rmJ Y\rangle$. Thus $\rmd\eta(e_1,e_2)=-2\tau\nu$.
\end{proof}

\begin{lemma} \label{diffT}
We have 
$$\rmd T^1=\sum_\gamma T^\gamma\omega^\gamma_1+\tau T^3\omega^2,$$
$$\rmd T^2=\sum_\gamma T^\gamma\omega^\gamma_2-\tau T^3\omega^1,$$
$$\rmd T^3=\sum_\gamma T^\gamma\omega^\gamma_3-\tau T^1\omega^2
+\tau T^2\omega^1.$$
\end{lemma}

\begin{proof}
The first two identities are a consequence of condition 
\eqref{conditionT1} and the last one of condition \eqref{conditionT2}.
\end{proof}

\begin{lemma} \label{diffOmega}
We have
\begin{eqnarray*}
\rmd\Omega+\Omega\wedge\Omega & = & \left(\begin{array}{ccc}
0 & \tau^2 & 0 \\
-\tau^2 & 0 & 0 \\
0 & 0 & 0
\end{array}\right)\omega^1\wedge\omega^2 \\
& & +(\kappa-4\tau^2)T^3\left(\begin{array}{ccc}
0 & T^3 & -T^2 \\
-T^3 & 0 & T^1 \\
T^2 & -T^1 & 0
\end{array}\right)\omega^1\wedge\omega^2.
\end{eqnarray*}
\end{lemma}

\begin{proof}
We set $\Psi=\rmd\Omega+\Omega\wedge\Omega$ and
$R^i_{klj}=\langle\rmR(e_k,e_l)e_j,e_i\rangle$. By proposition
\ref{differentiation} we have
$$\Psi^i_j=-\frac12\sum_k\sum_lR^i_{klj}\omega^k\wedge\omega^l
+\omega^i_3\wedge\omega^3_j,$$
and by the Gauss equation \eqref{gaussE} we have
$R^i_{klj}=\bar R^i_{klj}+\omega^3_j\wedge \omega^3_i(e_k,e_l)$ with
$$\bar R^i_{klj}=
(\kappa-3\tau^2)(\delta^k_j\delta^l_i-\delta^l_j\delta^k_i)
+(\kappa-4\tau^2)(T^lT^j\delta^k_i+T^kT^i\delta^l_j-T^lT^i\delta^k_j
-T^kT^j\delta^l_i).$$
Thus we get
$$\Psi^i_j=(\kappa-3\tau^2)\omega^i\wedge\omega^j
+(\kappa-4\tau^2)(T^i\omega^j-T^j\omega^i)\wedge\eta.$$

In the same way, by proposition \ref{differentiation} we have
$$\Psi^3_j=\frac12\sum_k\sum_l
\langle\nabla_{e_k}\rmS e_l-\nabla_{e_l}\rmS e_k-\rmS[e_k,e_l],e_j\rangle
\omega^k\wedge\omega^l,$$ 
and by the Codazzi equation \eqref{codazziE} we have
$$\langle\nabla_{e_k}\rmS e_l-\nabla_{e_l}\rmS e_k
-\rmS[e_k,e_l],e_j\rangle=
(\kappa-4\tau^2)T^3(T^l\delta^k_j-T^k\delta^l_j).$$
Thus we get 
$$\Psi^3_j=(\kappa-4\tau^2)T^3\omega^j\wedge\eta.$$

Hence we have
\begin{eqnarray*}
\Psi & = & (\kappa-3\tau^2)\left(\begin{array}{ccc}
0 & 1 & 0 \\
-1 & 0 & 0 \\
0 & 0 & 0
\end{array}\right)\omega^1\wedge\omega^2 \\
& & +(\kappa-4\tau^2)\left(\begin{array}{ccc}
0 & -T^2 & -T^3 \\
T^2 & 0 & 0 \\
T^3 & 0 & 0
\end{array}\right)\omega^1\wedge\eta \\
& & +(\kappa-4\tau^2)\left(\begin{array}{ccc}
0 & T^1 & -0 \\
-T^1 & 0 & -T^3 \\
0 & T^3 & 0
\end{array}\right)\omega^2\wedge\eta.
\end{eqnarray*}
We conclude using that $\omega^1\wedge\eta=T^2\omega^1\wedge\omega^2$,
$\omega^2\wedge\eta=-T^1\omega^1\wedge\omega^2$ and
$(T^1)^2+(T^2)^2+(T^3)^2=1$.
\end{proof}

\begin{lemma} \label{exprL}
We have 
\begin{eqnarray*}
L(Z) & = & (2\tau-\sigma)\left(
\begin{array}{ccc}
0 & -T^3 & T^2 \\
T^3 & 0 & -T^1 \\
-T^2 & T^1 & 0
\end{array}\right)\eta \\
& & +\left(
\begin{array}{ccc}
0 & 0 & 0 \\
0 & 0 & \tau \\
0 & -\tau & 0
\end{array}\right)\omega^1+\left(
\begin{array}{ccc}
0 & 0 & -\tau \\
0 & 0 & 0 \\
\tau & 0 & 0
\end{array}\right)\omega^2.
\end{eqnarray*}
\end{lemma}

\begin{proof}
We compute that
\begin{eqnarray*}
L(Z)^\alpha_\beta & = & 
\sum_k\left(\sum_\gamma\sum_\delta\sum_\varepsilon
Z^\varepsilon_\alpha Z^\gamma_k Z^\delta_\beta
\bar\Gamma^\delta_{\gamma\varepsilon}\right)\omega^k \\
& = & \sum_k(\tau(Z^2_\alpha Z^1_k Z^3_\beta
+Z^3_\alpha Z^2_k Z^1_\beta-Z^1_\alpha Z^2_k Z^3_\beta
-Z^3_\alpha Z^1_k Z^2_\beta) \\
& & \quad
+(\tau-\sigma)(Z^2_\alpha Z^3_k Z^1_\beta-Z^1_\alpha Z^3_k Z^2_\beta)
)\omega^k \\
& = & \sum_k(
\tau T^\beta(Z^1_k Z^2_\alpha-Z^1_\alpha Z^2_k)
+\tau T^\alpha(Z^1_\beta Z^2_k-Z^1_k Z^2_\beta) \\
& & \quad+(\tau-\sigma)T^k(Z^1_\beta Z^2_\alpha-Z^1_\alpha Z^2_\beta)
)\omega^k.
\end{eqnarray*}

Moreover the matrix $Z$ lies in $\mathrm{SO}_3(\R)$, so it is equal to
its comatrix. Using this fact we compute that
$$L(Z)^1_2=-(2\tau-\sigma)T^3(T^1\omega^1+T^2\omega^2),$$
$$L(Z)^1_3=(2\tau-\sigma)T^1T^2\omega^1+(2\tau-\sigma)(T^2)^2\omega^2
-\tau\omega^2,$$
$$L(Z)^2_3=-(2\tau-\sigma)(T^1)^2\omega^1-(2\tau-\sigma)T^1T^2\omega^2
+\tau\omega^1,$$
which proves the lemma.
\end{proof}

\begin{lemma} \label{LwedgeL}
We have 
\begin{eqnarray*}
L\wedge L & = & \tau(2\tau-\sigma)T^3\left(\begin{array}{ccc}
0 & -T^3 & T^2 \\
T^3 & 0 & -T^1 \\
-T^2 & T^1 & 0
\end{array}\right)\omega^1\wedge\omega^2 \\
& & +\tau(\tau-\sigma)\left(\begin{array}{ccc}
0 & 1 & 0 \\
-1 & 0 & 0 \\
0 & 0 & 0
\end{array}\right)\omega^1\wedge\omega^2.
\end{eqnarray*}
\end{lemma}

\begin{proof}
We compute that
\begin{eqnarray*}
L\wedge L & = & \tau(2\tau-\sigma)\left(\begin{array}{ccc}
0 & T^1 & 0 \\
-T^1 & 0 & -T^3 \\
0 & T^3 & 0
\end{array}\right)\eta\wedge\omega^2 \\
& & +\tau(2\tau-\sigma)\left(\begin{array}{ccc}
0 & -T^2 & -T^3 \\
T^2 & 0 & 0 \\
T^3 & 0 & 0
\end{array}\right)\eta\wedge\omega^1 \\
& & +\tau^2\left(\begin{array}{ccc}
0 & -1 & 0 \\
1 & 0 & 0 \\
0 & 0 & 0
\end{array}\right)\omega^1\wedge\omega^2.
\end{eqnarray*}
We conclude using that $(T^1)^2+(T^2)^2+(T^3)^2=1$.
\end{proof}

\begin{lemma} \label{LwedgeOmega}
We have 
\begin{eqnarray*}
L\wedge\Omega+\Omega\wedge L & = &
(2\tau-\sigma)\eta\wedge\left(\begin{array}{ccc}
0 & -\rmd T^3 & \rmd T^2 \\
\rmd T^3 & 0 & -\rmd T^1 \\
-\rmd T^2 & \rmd T^1 & 0
\end{array}\right) \\
& & +\tau(2\tau-\sigma)T^3\left(\begin{array}{ccc}
0 & T^3 & -T^2 \\
-T^3 & 0 & T^1 \\
T^2 & -T^1 & 0
\end{array}\right)\omega^1\wedge\omega^2 \\
& & +\tau(2\tau-\sigma)\left(\begin{array}{ccc}
0 & -1 & 0 \\
1 & 0 & 0 \\
0 & 0 & 0
\end{array}\right)\omega^1\wedge\omega^2 \\
& & +\tau\left(\begin{array}{ccc}
0 & 0 & 0 \\
0 & 0 & -1 \\
0 & 1 & 0
\end{array}\right)\rmd\omega^1
+\tau\left(\begin{array}{ccc}
0 & 0 & 1 \\
0 & 0 & 0 \\
-1 & 0 & 0
\end{array}\right)\rmd\omega^2.
\end{eqnarray*}
\end{lemma}

\begin{proof}
We compute that
\begin{eqnarray*}
L\wedge\Omega+\Omega\wedge L & = &
(2\tau-\sigma)\eta\wedge M \\
& & +\tau\omega^2\wedge\left(\begin{array}{ccc}
0 & -\omega^3_2 & 0 \\
-\omega^2_3 & 0 & \omega^2_1 \\
0 & \omega^1_2 & 0
\end{array}\right) \\
& & +\tau\omega^1\wedge\left(\begin{array}{ccc}
0 & \omega^1_3 & -\omega^1_2 \\
\omega^3_1 & 0 & 0 \\
-\omega^2_1 & 0 & 0
\end{array}\right)
\end{eqnarray*}
with $$M=\left(\begin{array}{ccc}
0 & T^2\omega^3_2-T^1\omega^1_3 & -T^3\omega^2_3+T^1\omega^1_2 \\
-T^1\omega^3_1+T^2\omega^2_3 & 0 & T^3\omega^1_3-T^2\omega^2_1 \\
T^1\omega^2_1-T^3\omega^3_2 & -T^2\omega^1_2+T^3\omega^3_1 & 0
\end{array}\right).$$

We conclude using lemma \ref{diffT}, formulas \eqref{diffomega1} and
\eqref{diffomega2}, and the fact that $(T^1)^2+(T^2)^2+(T^3)^2=1$.
\end{proof}

For $y\in\cV$, let $\cZ(y)$ be the set of matrices
$Z\in\mathrm{SO}_3(\R)$ such 
that the coefficients of the last line of $Z$ are the $T^\beta(y)$.
It is diffeomorphic to the circle $\s^1$.


We now prove the following proposition.

\begin{prop} \label{matrixA}
Assume that the compatibility equations for $\E$ are
satisfied. Let $y_0\in\cV$ and $A_0\in\cZ(y_0)$. Then there exist a
neighbourhood $U_1$ of
$y_0$ in $\cV$ and a unique map $A:U_1\to\mathrm{SO}_3(\R)$ such that
$$A^{-1}\rmd A=\Omega,$$
$$\forall y\in U_1,\quad A(y)\in\cZ(y),$$
$$A(y_0)=A_0.$$
\end{prop}

\begin{proof}
Let $U$ be a coordinate neighbourhood in $\cV$. The set 
$$\cF=\{(y,Z)\in U\times\mathrm{SO}_3(\R);Z\in\cZ(y)\}$$ is a
manifold of dimension $3$, and
$$\rmT_{(y,Z)}\cF=\{(u,\zeta)\in\rmT_yU\oplus\rmT_Z\mathrm{SO}_3(\R);
\zeta^3_\beta=(\rmd T^\beta)_y(u)\}.$$

Let $Z$ denote the projection
$U\times\mathrm{SO}_3(\R)
\to\mathrm{SO}_3(\R)\subset\cM_3(\R)$.
We consider on $\cF$ the 
following matrix of $1$-forms:
$$\Theta=Z^{-1}\rmd Z-\Omega-L(Z)$$
where $L(Z)$ is defined in lemma \ref{diffA}, 
namely for $(y,Z)\in\cF$ we have
$$\Theta_{(y,Z)}:\rmT_{(y,Z)}\cF\to\cM_3(\R),$$
$$\Theta_{(y,Z)}(u,\zeta)=Z^{-1}\zeta-\Omega_y(u)-L(Z)(u).$$

We claim that, for each $(y,Z)\in\cF$, the space
$$\cD(y,Z)=\ker\Theta_{(y,Z)}$$ has dimension $2$.
We first notice that
the matrix $\Theta$ belongs to $\mathfrak{so}_3(\R)$ since
$\Omega$, $L(Z)$ and $Z^{-1}\rmd Z$ do. Moreover we have
$$(Z\Theta)^3_\beta
=\rmd Z^3_\beta-\sum_\gamma Z^3_\gamma\omega^\gamma_\beta
-\sum_\gamma Z^3_\gamma L(Z)^\gamma_\beta
=\rmd T^\beta-\sum_\gamma T^\gamma\omega^\gamma_\beta
-\sum_\gamma T^\gamma L(Z)^\gamma_\beta.$$
Using lemmas \ref{diffT} and \ref{exprL} we compute that
$$(Z\Theta)^3_\beta=0.$$
Thus the values of
$\Theta_{(y,Z)}$ lie in the space
$$\cH=\{H\in\mathfrak{so}_3(\R);(ZH)^3_\beta=0\},$$
which has dimension $1$ (indeed, the map 
$F:\mathrm{SO}_3(\R)\to\s^2,
Z\mapsto(Z^3_\beta)_\beta$ is a submersion, and we
have $H\in\cH$ if and 
only if $ZH\in\ker(\rmd F)_Z$).
Moreover, the space $\rmT_{(y,Z)}\cF$ contains the subspace
$\{(0,ZH);H\in\cH\}$, and the restriction of $\Theta_{(y,Z)}$ on this
subspace is the map $(0,ZH)\mapsto H$. Thus
$\Theta_{(y,Z)}$ is onto 
$\cH$, and consequently the linear map $\Theta_{(y,Z)}$ has rank
$1$. This finishes proving the claim.

We now prove that the distribution $\cD$ is involutive. We first compute
that 
\begin{eqnarray*}
\rmd\Theta & = & -Z^{-1}\rmd Z\wedge Z^{-1}\rmd Z-\rmd\Omega-\rmd L \\
& = & -(\Theta+\Omega+L)\wedge(\Theta+\Omega+L)-\rmd\Omega-\rmd L \\
& = & -\Theta\wedge\Theta-\Theta\wedge\Omega
-\Omega\wedge\Theta-\Omega\wedge L-L\wedge\Omega \\
& & -\Omega\wedge\Omega-\rmd\Omega-L\wedge L-\rmd L.
\end{eqnarray*}
Using lemmas \ref{diffeta}, \ref{diffOmega}, \ref{LwedgeL},
\ref{LwedgeOmega} and the relation $\sigma=\frac\kappa{2\tau}$, we
obtain
$$\rmd\Theta=-\Theta\wedge\Theta-\Theta\wedge\Omega-\Omega\wedge\Theta.$$
From this formula we deduce that if $\xi_1,\xi_2\in\cD$, then
$\rmd\Theta(\xi_1,\xi_2)=0$, and so
$\Theta([\xi_1,\xi_2])=\xi_1\cdot\Theta(\xi_2)
-\xi_2\cdot\Theta(\xi_1)-\rmd\Theta(\xi_1,\xi_2)=0$, i.e.,
$[\xi_1,\xi_2]\in\cD$. Thus the distribution $\cD$ is involutive, and
so, by the theorem of Frobenius, it is integrable.

Let $\cA$ be the integral manifold through $(y_0,A_0)$.
If $\zeta\in\rmT_{A_0}\mathrm{SO}_3(\R)$ is such that
$(0,\zeta)\in\rmT_{(y_0,A_0)}\cA=\cD(y_0,A_0)$, 
then we have $0=\Theta_{(y_0,A_0)}(0,\zeta)=A_0^{-1}\zeta$. This
proves that
$$\rmT_{(y_0,A_0)}\cA\cap
\left(\{0\}\times\rmT_{A_0}\mathrm{SO}_3(\R)\right)=\{0\}.$$
Thus the manifold
$\cA$ is locally the graph of a function
$A:U_1\to\mathrm{SO}_3(\R)$ where $U_1$ is a neighbourhood of
$y_0$ in $U$. By construction, this map satisfies the properties of
proposition \ref{matrixA} and is unique. 
\end{proof}

\begin{prop} \label{functionf}
Let $x_0\in\E$ (without loss of generality we can assume that 
$x_0\in\E'$). There exist a neighbourhood $U_2$ of $y_0$ contained in
$U_1$ and a unique function $f:U_2\to\E'$ such that
$$\rmd f=(B\circ f)A\omega,$$
$$f(y_0)=x_0,$$
where $\omega$ is the column $(\omega^1,\omega^2,0)$ and, for $x\in\E'$,
$B(x)\in\cM_3(\R)$ is the matrix of the coordinates of the frame
$(E_\alpha(x))$ in the frame 
$(\partial_{x^\alpha})$.
\end{prop}

\begin{proof}
We consider on $U_1\times\E'$ the following matrix of $1$-forms:
$$\Lambda=B^{-1}\rmd x-A\omega,$$
namely, for $q\in U_1$ and $x\in\E'$ we have
$$\Lambda_{(q,x)}:\rmT_qU_1\oplus\rmT_x\E\to\cM_{3,1}(\R),$$
$$\Lambda_{(q,x)}(u,v)=B(x)^{-1}v-A(q)\omega_q(u).$$

We first notice that for all $(q,x)\in U_1\times\E'$ the linear
map $\Lambda_{(q,x)}$ is onto $\cM_{3,1}(\R)$. Consequently the space
$$\cE(q,x)=\ker\Lambda_{(q,x)}$$ has dimension $2$. We will prove that this
distribution $\cE$ is integrable.

We have
$$\rmd\Lambda=-B^{-1}\rmd BB^{-1}\wedge\rmd x
-\rmd A\wedge\omega-A\wedge\rmd\omega.$$
By equations \eqref{diffomega1} and \eqref{diffomega2} we have
$\rmd\omega=-\Omega\wedge\omega$; and by proposition \ref{matrixA} we
have $\rmd A=A\Omega+AL(A)$. Thus we get
$$\rmd\Lambda=-B^{-1}\rmd B\wedge\Lambda-B^{-1}\rmd B\wedge A\omega
-AL(A)\wedge\omega.$$

Using lemma \ref{exprL} we compute that
$$L(A)\wedge\omega=-(2\tau-\sigma)T^3\left(\begin{array}{c}
T^1 \\
T^2 \\
T^3
\end{array}\right)\omega^1\wedge\omega^2
-\left(\begin{array}{c}
0 \\
0 \\
\sigma
\end{array}\right)\omega^1\wedge\omega^2,$$
and thus, using the fact that $A^3_\beta=T^\beta$ and 
$A=\mathrm{com}A$, we get
$$AL(A)\wedge\omega=\left(\begin{array}{c}
-\sigma A^1_3 \\
-\sigma A^2_3 \\
-2\tau T^3
\end{array}\right)\omega^1\wedge\omega^2.$$

%

We will use the notation $(x,y,x)$ instead of $(x^1,x^2,x^3)$ for the
coordinates in $\E$ and we will use the local models described in
sections \ref{bergerspheres}, \ref{heisenberg} and \ref{psl}. Using
formulas \eqref{canonicalbergerspheres}, 
\eqref{canonicalheisenberg} and \eqref{canonicalpsl}, we get that the
matrix $B$ is
$$B=\left(\begin{array}{ccc}
\lambda^{-1}\cos(\sigma z) & -\lambda^{-1}\sin(\sigma z) & 0 \\
\lambda^{-1}\sin(\sigma z) & \lambda^{-1}\cos(\sigma z) & 0 \\
\tau(x\sin\sigma z-y\cos\sigma z) & \tau(x\cos\sigma z+y\sin\sigma z)
& 1
\end{array}\right),$$
with $$\lambda=\frac1{1+\frac\kappa4(x^2+y^2)}.$$

We will write $$A\omega=\left(\begin{array}{c}
\alpha^1 \\
\alpha^2 \\
\eta
\end{array}\right)$$ with
$$\alpha^j=A^j_1\omega^1+A^j_2\omega^2.$$
Then we have
$$\Lambda=B^{-1}\rmd X-A\omega=\left(\begin{array}{c}
\lambda(\cos(\sigma z)\rmd x+\sin(\sigma z)\rmd y)-\alpha^1 \\
\lambda(-\sin(\sigma z)\rmd x+\cos(\sigma z)\rmd y)-\alpha^2 \\
\tau\lambda(y\rmd x-x\rmd y)+\rmd z-\eta
\end{array}\right).$$
We also compute that
$$B^{-1}\rmd B=\left(\begin{array}{ccc}
\frac\kappa2\lambda(x\rmd x+y\rmd y) & -\sigma\rmd z & 0 \\
\sigma\rmd z & \frac\kappa2\lambda(x\rmd x+y\rmd y) & 0 \\
a & b & 0
\end{array}\right)$$
with $$a=\frac{\tau\kappa}2\lambda(y\cos(\sigma z)-x\sin(\sigma z))
(x\rmd x+y\rmd y)+\tau(\sin(\sigma z)\rmd x-\cos(\sigma z)\rmd y),$$
$$b=-\frac{\tau\kappa}2\lambda(x\cos(\sigma z)+y\sin(\sigma z))
(x\rmd x+y\rmd y)+\tau(\cos(\sigma z)\rmd x+\sin(\sigma z)\rmd y).$$
Thus we have
\begin{eqnarray*}
B^{-1}\rmd B\wedge A\omega+AL(A)\wedge\omega & = &
\left(\begin{array}{c}
\frac\kappa2\lambda(x\rmd x+y\rmd y)\wedge\alpha^1
-\sigma\rmd z\wedge\alpha^2 \\
\sigma\rmd z\wedge\alpha^1
+\frac\kappa2\lambda(x\rmd x+y\rmd y)\wedge\alpha^2 \\
a\wedge\alpha^1+b\wedge\alpha^2
\end{array}\right) \\
& & +\left(\begin{array}{c}
-\sigma A^1_3 \\
-\sigma A^2_3 \\
-2\tau T^3
\end{array}\right)\omega^1\wedge\omega^2.
\end{eqnarray*}

Using the above expression for $\Lambda$ we get
$$\lambda\rmd x=\cos(\sigma z)\Lambda^1-\sin(\sigma z)\Lambda^2
+\cos(\sigma z)\alpha^1-\sin(\sigma z)\alpha^2,$$
$$\lambda\rmd y=\sin(\sigma z)\Lambda^1+\cos(\sigma z)\Lambda^2
+\sin(\sigma z)\alpha^1+\sin(\sigma z)\alpha^2,$$
$$\rmd z=\Lambda^3+\eta-\tau\lambda(y\rmd x-x\rmd y).$$

The term in the first line of the matrix 
$B^{-1}\rmd B\wedge A\omega+AL(A)$ is 
\begin{eqnarray*}
\frac\kappa2(y\cos(\sigma z)-x\sin(\sigma z))\alpha^2\wedge\alpha^1
+\sigma\tau(y\cos(\sigma z)-x\sin(\sigma z))\alpha^1\wedge\alpha^2 \\
\quad-\sigma\eta\wedge\alpha^2-\sigma A^1_3\omega^1\wedge\omega^2
+\chi^1
\end{eqnarray*}
where $\chi^1$ is a linear combination of the $\Lambda^\alpha$ (the
coefficients being $1$-forms).
Since $\sigma=\frac{\kappa}{2\tau}$, the first two terms in this expression
cancel. Moreover we have $\eta\wedge\alpha^2
=(A^3_1A^2_2-A^3_2A^2_1)\omega^1\wedge\omega^2
=-A^1_3\omega^1\wedge\omega^2$, hence the term in the first line of
the matrix $B^{-1}\rmd B\wedge A\omega+AL(A)$ is $\chi^1$. In the same
way, the term in the second line of
the matrix $B^{-1}\rmd B\wedge A\omega+AL(A)$ is 
a linear combination of the $\Lambda^\alpha$ which will be
denoted by $\chi^2$. Finally we compute that the term in the
third line of the matrix $B^{-1}\rmd B\wedge A\omega+AL(A)$ is 
$$\left(\frac{2\tau}{\lambda}-\frac{\tau\kappa}2(x^2+y^2)\right)
\alpha^1\wedge\alpha^2-2\tau T^3\omega^1\wedge\omega^2+\chi^3$$
where $\chi^1$ is a linear combination of the $\Lambda^\alpha$.
Since $\lambda^{-1}=1+\frac\kappa4(x^2+y^2)$ and 
$\alpha^1\wedge\alpha^2=(A^1_1A^2_2-A^1_2A^2_1)\omega^1\wedge\omega^2
=T^3\omega^1\wedge\omega^2$, this term is simply $\chi^3$.
We conclude that
$$B^{-1}\rmd B\wedge A\omega+AL(A)=\chi$$
where $\chi$ is a matrix of $2$-forms which are linear combinations
of the coefficients of $\Lambda$. Finally we have
$$\rmd\Lambda=-B^{-1}\rmd B\wedge\Lambda-\chi.$$
From this formula we deduce that if $\xi_1,\xi_2\in\cE$, then
$\rmd\Lambda(\xi_1,\xi_2)=0$, and so
$[\xi_1,\xi_2]\in\cE$. Thus the distribution $\cE$ is involutive, and
so, by the theorem of Frobenius, it is integrable.

Let $\cA$ be the integral manifold through $(y_0,x_0)$.
If $v\in\rmT_{x_0}\E$ is such that
$(0,v)\in\rmT_{(y_0,x_0)}\cA=\cD(y_0,x_0)$, 
then we have $0=\Lambda_{(y_0,x_0)}(0,v)=B(x_0)^{-1}v$. This
proves that
$$\rmT_{(y_0,x_0)}\cA\cap
\left(\{0\}\times\rmT_{x_0}\E\right)=\{0\}.$$
Thus the manifold
$\cA$ is locally the graph of a function
$A:U_2\to\E'$ where $U_2$ is a neighbourhood of
$y_0$ in $U_1$. By construction, this map satisfies the properties of
proposition \ref{matrixA} and is unique.



\end{proof}

We now prove the theorem.

\begin{proof}[Proof of theorem \ref{isometry}]
Let $y_0\in\cV$, $A_0\in\cZ(y_0)$ and $x_0\in\E'$.
We consider on $\cV$ a local orthonormal frame $(e_1,e_2)$ in
the neighbourhood of $y_0$ and we keep the same notations. Then by
propositions \ref{matrixA} and \ref{functionf} there exists a unique map 
$A:U_2\to\mathrm{SO}^3(\R)$ such that
$$A^{-1}\rmd A=\Omega+L(A),$$
$$\forall y\in U_1,\quad A(y)\in\cZ(y),$$
$$A(y_0)=A_0,$$
and a unique map $f:U_2\to\E'$ such that
$$\rmd f=(B\circ f)A\omega,$$
$$f(y_0)=x_0,$$
 where $U_2$ is a neighbourhood of $y_0$, which we can
assume simply connected. We will check that $f$ has the properties
required in the theorem on $U_2$.

We have $\rmd f^\alpha(e_k)=(B(f)A)^\alpha_k$, so in the frame 
$(\partial_{x^\alpha})$ the vector
$\rmd f(e_k)$ is given by the column $k$ of the matrix $BA$, which is
invertible. Hence $\rmd f$ has rank $2$, and thus $f$ is an immersion.
Moreover, in the frame $(E_\alpha)$ the vector
$\rmd f(e_k)$ is given by the column $k$ of the matrix $A$, which is
orthogonal, and thus we have
$\langle\rmd f(e_p),\rmd f(e_q)\rangle
=\delta^p_q$, which means that $f$ is an isometry.

The columns of $A(y)$ form a direct orthonormal frame of
$\E$. The first and second columns form a direct orthonormal frame of
$\rmT_{f(y)}f(\cV)$ 
Thus the third column gives, in the frame $(E_\alpha)$,
the unit normal $N(f(y))$ to
$f(\cV)$ in $\E$ at the point $f(y)$.

We set $X_j=\rmd f(e_j)$. Then we have
\begin{eqnarray*}
\rmd A^\alpha_j(e_k) & = & \langle\bar\nabla_{X_k}X_j,E_\alpha\rangle
+\langle X_j,\bar\nabla_{X_k}E_\alpha\rangle \\
& = & \langle\bar\nabla_{X_k}X_j,E_\alpha\rangle
+\sum_\gamma\sum_\delta A^\gamma_kA^\delta_j
\bar\Gamma^\delta_{\gamma\alpha} \\
& = & \langle\bar\nabla_{X_k}X_j,E_\alpha\rangle
+(AL(A))^\alpha_j(e_k),
\end{eqnarray*}
so
\begin{eqnarray*}
\langle\bar\nabla_{X_k}X_j,N\rangle & = &
\sum_\alpha\langle\bar\nabla_{X_k}X_j,E_\alpha\rangle A^\alpha_3
=\sum_\alpha A^\alpha_3(\rmd A-AL(A))^\alpha_j(e_k) \\
& = & \sum_\alpha A^\alpha_3(A\Omega)^\alpha_j(e_k)
=\sum_\alpha\sum_\gamma A^\alpha_\gamma
A^\alpha_3\omega^\gamma_j(e_k) \\
& = & \omega^3_j(e_k)=\langle\rmS e_k,e_j\rangle.
\end{eqnarray*}
This means that the shape operator of $f(\cV)$ in $\E$ is 
$\rmd f\circ\rmS\circ\rmd f^{-1}$.

Finally, the coefficients of the vertical vector
$\xi=E_3$ in 
the orthonormal frame $(X_1,X_2,N)$ are given by the last
line of $A$. Since $A(y)\in\cZ(y)$ for all $y\in U_2$ we get
$$\xi=\sum_j T^jX_j+T^3N
=\rmd f(T)+\nu N.$$

We now prove that the local immersion is unique up to a global
isometry of $\E$ preserving $\xi$ (and also, consequently, the 
orientation of the base of the fibration). Let $\tilde f:U_3\to\E$ be
another immersion satisfying the conclusion of the theorem, where
$U_3$ is a simply connected neighbourhood of $y_0$ included in $U_2$,
let $(\tilde X_\beta)$ be the associated frame (i.e., $\tilde
X_j=\rmd\tilde f(e_j)$ and $\tilde X_3$ is the normal of 
$\tilde f(\cV)$) and let $\tilde A$ the
matrix of the coordinates of the frame $(\tilde X_\beta)$ in the
frame $(E_\alpha)$. Up to an isometry of $\E$ (which is necessarily
direct), we can
assume that $f(y_0)=\tilde f(y_0)$ and that the frames
$(X_\beta(y_0))$ and $(\tilde X_\beta(y_0))$ coincide, i.e.,
$A(y_0)=\tilde A(y_0)$. We notice that this isometry necessarily fixes
$\xi$ since the $T^\alpha$ are the same for $x$
and $\tilde x$. The matrices $A$ and $\tilde A$ satisfy
$A^{-1}\rmd A=\Omega+L(A)$ and 
$\tilde A^{-1}\rmd\tilde A=\Omega+L(\tilde A)$ (see
section \ref{hypersurfaces}),
$A(y),\tilde A(y)\in\cZ(y)$ and $A(y_0)=\tilde A(y_0)$,
thus by the uniqueness of the solution of the equation in proposition
\ref{matrixA} we get $A(y)=\tilde A(y)$. We conclude similarly that
$f=\tilde f$ on $U_3$.

The proof that this local immersion $f$ can be extended to the whole $\cV$
(since $\cV$ is simply connected) is exactly the same as the proof
of the corresponding statement in theorem 3.3 in \cite{codazzi} (it is
a standard argument).
\end{proof}

%

\begin{rem} \label{changeofsigns}
If $(\rmd s^2,\rmS,T,\nu)$ satisfies the compatibilty equations and
correspond to an immerion $f:\Sigma\to\E$, then 
$(\rmd s^2,\rmS,-T,-\nu)$ also satisfies the compatibilty equations and
corresponds to the immersion $\sigma\circ f$ where $\sigma$ is an
isometry of $\E$ 
reversing the orientations of both the fibers and the base of the
fibration.
\end{rem}

\section{Constant mean curvature surfaces in $3$-dimensional homogeneous manifolds}

In this section we will give an application of theorem \ref{isometry}
to constant mean curvature surfaces (CMC) in $3$-dimensional homogeneous manifolds with
$4$-dimensional isometry group. Abresch and Rosenberg proved that there
exists a holomorphic quadratic differential for CMC surfaces in
$\s^2\times\R$ and $\h^2\times\R$, generalizing the Hopf differential for
CMC surfaces in $3$-dimensional space forms (\cite{abresch}). Since the
Hopf differential is a very useful tool for CMC surfaces, this
motivated many works on CMC surfaces in $\s^2\times\R$ and $\h^2\times\R$.
Recently, Abresch announced the existence of a holomorphic quadratic
differential for CMC surfaces in all $3$-dimensional homogeneous manifolds with
$4$-dimensional isometry group (\cite{abreschsurvey}). 
This indicates that the theory of CMC
surfaces in these manifolds may be particularily interesting.

We will consider constant mean curvature immersions of
oriented surfaces. Consequently the 
mean curvature will be defined
with a sign: it will be positive
if the mean curvature vector induces the same orientation as the initial
orientation, and it will be negative
if the mean curvature vector induces the opposite orientation.

We will denote by $\rmI$ and $\rmJ$ the identity and the rotation of
angle $\frac\pi2$ on the tangent bundle of a surface.

\subsection{A generalized Lawson correspondence}

It is well known that there exists an isometric correspondence 
between certain simply connected CMC surfaces in space-forms:
more precisely, every simply connected CMC $H_1$ surface in
$\M^3(K_1)$ is isometric to a simply connected CMC $H_2$ surface in
$\M^3(K_2)$ with $K_1-K_2=H_2^2-H_1^2$, and the shape operators of these
two surfaces differ by $(H_2-H_1)\rmI$.
Two such surfaces are called
cousin surfaces. This correspondence is often called the Lawson
correspondence. In particular, any simply connected minimal surface
in $\s^3$ is isometric to a CMC $1$ surface in $\R^3$, and any
minimal surface
in $\R^3$ is isometric to a CMC $1$ surface in $\h^3$.

The Lawson correspondence is a consequence of the Gauss and Codazzi
equations in the space-forms.

In this section we will use the compatibility equations for
homogeneous $3$-manifolds with $4$-dimensional isometry group
and theorem \ref{isometry} to prove the existence of an
isometric correspondence between certain simply connected
CMC surfaces in these $3$-manifolds. Hence this will be a generalisation
of the Lawson correspondence.

The technique will be to start with some data $(\rmd s^2,\rmS,T,\nu)$
on a surface satisfying the compatibility equations for some homogeneous
$3$-manifold and to modify them in order to get data
satisfying the compatibility equations for another homogeneous
$3$-manifold. An important fact is that the space of symmetric traceless
operators is globally invariant by rotation.
The easiest change is to keep $\rmd s^2$ and $\nu$, and to rotate $T$
and the traceless part of $\rmS$ by some fixed angles; the Codazzi
equation then implies that we need to take the same angle for $T$ and
the traceless part of $\rmS$.

\begin{prop} \label{correspondence}
Let $\E_1$ and $\E_2$ be two $3$-dimensional homogeneous manifolds
with $4$-dimensional isometry groups, of base curvatures $\kappa_1$
and $\kappa_2$ and bundle curvatures $\tau_1$ and $\tau_2$
respectively. Assume that
$$\kappa_1-4\tau_1^2=\kappa_2-4\tau_2^2.$$

Let $H_1$ and $H_2$ be two real numbers such that
$$\tau_1^2+H_1^2=\tau_2^2+H_2^2.$$

Let $\cV$ be a surface with a quadruple $(\rmd s^2,\rmS_1,T_1,\nu)$
satisfying the compatibility equations for $\E_1$ and such that 
$$\tr\rmS_1=2H_1.$$

Let $$\theta\in\R,$$ 
$$T_2=e^{\theta\rmJ}T_1,$$
$$\rmS_2=e^{\theta\rmJ}(\rmS_1-H_1\rmI)+H_2\rmI.$$
In particular $\rmS_2$ is symmetric and satisfies $$\tr\rmS_2=2H_2.$$

If the real number $\theta$ satisfies
\begin{equation} \label{phase}
\tau_2+iH_2=e^{i\theta}(\tau_1+iH_1),
\end{equation}
then the quadruple $(\rmd s^2,\rmS_2,T_2,\nu)$
satisfies the compatibility equations for $\E_2$.

Conversely, if the function $\nu$ is not identically zero
and if the quadruple $(\rmd s^2,\rmS_2,T_2,\nu)$
satisfies the compatibility equations for $\E_2$, then
\eqref{phase} holds.
\end{prop}

\begin{proof}
The fact that $\rmS_2$ is symmetric comes from the fact that the space
of symmetric traceless operators is invariant by a rotation.

We have $$\det(\rmS_k-H_k\rmI)=\det\rmS_k-H_k^2$$ for $k=1,2$,
and so $$\det\rmS_1=\det\rmS_2+H_1^2-H_2^2.$$
Let $K$ be the Gauss curvature of the metric $\rmd s^2$.
By the Gauss equation \eqref{gaussE} we have
\begin{eqnarray*}
K & = & \det\rmS_1+\tau_1^2+(\kappa_1-4\tau_1^2)\nu^2 \\
& = & \det\rmS_2+H_1^2-H_2^2+\tau_1^2+(\kappa_1-4\tau_1^2)\nu^2 \\
& = & \det\rmS_2+\tau_2^2+(\kappa_2-4\tau_2^2)\nu^2
\end{eqnarray*}
since $\kappa_1-4\tau_1^2=\kappa_2-4\tau_2^2$ and
$\tau_1^2+H_1^2=\tau_2^2+H_2^2$. Thus the quadruple
$(\rmd s^2,\rmS_2,T_2,\nu)$ satisfies the Gauss equation for
$\E_2$.

Since
$\rmJ$ commutes with $\nabla_X$ for all vector fields $X$, we have
$$\nabla_X\rmS_2 Y-\nabla_Y\rmS_2 X-\rmS_2[X,Y]=
e^{\theta\rmJ}(\nabla_X\rmS_1 Y-\nabla_Y\rmS_1 X-\rmS_1[X,Y]).$$
On the other hand, a computation done in the proof of proposition 4.1 in
\cite{codazzi} shows that
$$\langle Y,T_2\rangle X-\langle X,T_2\rangle Y=
e^{\theta\rmJ}(\langle Y,T_1\rangle X-\langle X,T_1\rangle Y).$$
Hence the Codazzi equation
for $\E_2$ is satisfied by $(\rmd s^2,\rmS_2,T_2,\nu)$.

To prove that the quadruple $(\rmd s^2,\rmS_2,T_2,\nu)$
satisfies the compatibility equations \eqref{conditionT1}
and \eqref{conditionT2} for $\E_2$, it suffices to prove that
\begin{equation} \label{rotationshape}
\rmS_2-\tau_2\rmJ=e^{\theta\rmJ}(\rmS_1-\tau_1\rmJ).
\end{equation}
Using the expression of $\rmS_2$, equation \eqref{rotationshape}
is equivalent to
\begin{equation} \label{rotationshape2}
H_2\rmI-\tau_2\rmJ=e^{\theta\rmJ}(H_1\rmI-\tau_1\rmJ).
\end{equation}
We notice that this is a purely algebraic condition: the shape
operators are not involved anymore.
We consider a local direct orthonormal frame and we will 
identify the operators with their matrix in this frame.
Then we have
$$\rmJ=\left(\begin{array}{cc}
0 & -1 \\
1 & 0
\end{array}\right).$$ Then equation \eqref{rotationshape2}
is equivalent to
$$\left\{\begin{array}{ccc}
H_2 & = & H_1\cos\theta+\tau_1\sin\theta, \\
\tau_2 & = & \tau_1\cos\theta-H_1\sin\theta.
\end{array}\right.,$$
i.e., it is equivalent to equation \eqref{phase}. This proves
the first assertion of the theorem.

Conversely, if $(\rmd s^2,\rmS_2,T_2,\nu)$ satisfies the compatibility
equations for $\E_2$, then the compatibility equations
\eqref{conditionT1} for $(\rmd s^2,\rmS_1,T_1,\nu)$ and
$(\rmd s^2,\rmS_2,T_2,\nu)$ imply that
\eqref{rotationshape} holds at every point where $\nu\neq 0$.
If there exists a point where $\nu\neq 0$, this implies that
\eqref{phase} holds. 
\end{proof}

\begin{thm} \label{sisters}
Let $\E_1$ and $\E_2$ be two $3$-dimensional homogeneous manifolds
with $4$-dimensional isometry groups, of base curvatures $\kappa_1$
and $\kappa_2$ and bundle curvatures $\tau_1$ and $\tau_2$
respectively, and such that
$$\kappa_1-4\tau_1^2=\kappa_2-4\tau_2^2.$$
Let $\xi_1$ and $\xi_2$ be the vertical vector fields of $\E_1$ and
$\E_2$ respectively.

Let $\Sigma$ be a simply connected
Riemann surface and let $x_1:\Sigma\to\E_1$ be
a conformal constant mean curvature $H_1$ immersion with
$H_1^2\geqslant\tau_2^2-\tau_1^2$.
Let $N_1$ be the induced normal (compatible with the orientation
of $\Sigma$). Let
$\rmS_1$ be the symmetric operator on $\Sigma$ induced by the shape
operator of $x_1(\Sigma)$ associated to the normal $N_1$. Let $T_1$
be the vector field on $\Sigma$ such that $\rmd x_1(T_1)$ is the
projection of $\xi_1$ onto $\rmT(x_1(\Sigma))$. Let
$\nu=\langle N_1,\xi_1\rangle$. 

Let $H_2\in\R$ such that
$$\tau_1^2+H_1^2=\tau_2^2+H_2^2.$$

Let $\theta\in\R$ such that
$$\tau_2+iH_2=e^{i\theta}(\tau_1+iH_1).$$

Then there exists a conformal immersion
$x_2:\Sigma\to\E_2$ such that:
\begin{enumerate}
\item the metrics induced on $\Sigma$ by $x_1$ and $x_2$ are the same,
\item the symmetric operator on $\Sigma$ induced by the shape operator
of $x_2(\Sigma)$ is $e^{\theta\rmJ}(\rmS_1-H_1\rmI)+H_2\rmI$,
\item $\xi_2=\rmd x_2(e^{\theta\rmJ}T_1)+\nu N_2$ where 
$N_2$ is the unit normal to $x_2$.
\end{enumerate}

Moreover, this immersion $x_2$ is unique up to isometries of $\E_2$
preserving the orientations of both the fibers and the base of the
fibration, and it has constant mean curvature $H_2$.

The immersions $x_1$ and $x_2$ are called sister immersions. The 
number $\theta$ is called the phase of $(x_1,x_2)$.
\end{thm}

This means that there exists an isometric
correspondence between CMC $H_1$ simply connected surfaces
in $\E_1$ and CMC $H_2$ simply connected surfaces in $\E_2$.

\begin{proof}
Let $\rmd s^2$ be the metric on $\Sigma$ induced by $x_1$. Then $(\rmd
s^2,\rmS_1,T_1,\nu)$ satisfies the compatibility equations for
$\E_1$. Thus, by proposition \ref{correspondence}, the quadruple
$(\rmd s^2,\rmS_2,e^{\theta\rmJ}T_1,\nu)$
with $\rmS_2=e^{\theta\rmJ}(\rmS_1-H_1\rmI)+H_2\rmI$ also does. Thus by
theorem \ref{isometry} there exists an immersion $x_2$
satisfying properties 1, 2, and 3, and this immersion is unique up to
isometries of $\E_2$
preserving the orientations of both the fibers and the base of the
fibration. 
Moreover, we have
$\tr\rmS_2=2H_2$, i.e., the immersion
$x_2$ has mean curvature $H_2$. 
\end{proof}

\begin{figure}[htbp]
\begin{center}
\input{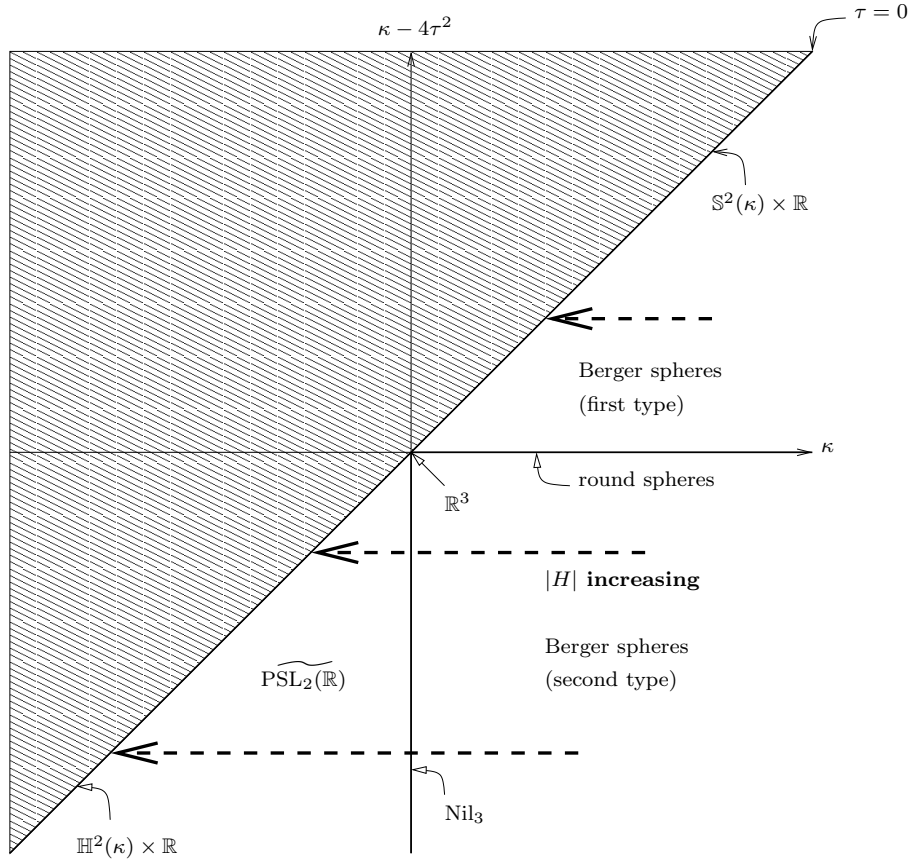}
\caption{The correspondence of the sister surfaces}
\label{figuresisters}
\end{center}
\end{figure}

Figure \ref{figuresisters} helps visualizing which classes of CMC surfaces
are related by the sister surface correspondence. We start from a CMC
surface in some homogeneous $3$-manifold. 
Then we can go horizontally on
the graph. We can go to the left until reaching a manifold with $\tau=0$;
in this case the absolute mean curvature
$|H|$ increases. We can go to the right until reaching
$H=0$; in this case $|H|$ decreases.

A particularily interesting case is when $\E_1$ is the Heisenberg space
$\nil$ with its standard metric ($\kappa_1=0$, $\tau_1=\frac12$) and
$\E_2=\h^2\times\R$ ($\kappa_2=-1$, $\tau_2=0$). Then CMC $H_1$
surfaces in $\nil$ correspond isometrically
to CMC $H_2$ surfaces in $\h^2\times\R$ with
$H_2^2=H_1^2+\frac14$. In particular we have the following corollary.

\begin{cor} \label{sistersheisenberg}
There exists an isometric correspondence with phase $\theta=\frac\pi2$
between simply connected
minimal surfaces in the Heisenberg space $\nil$ and
simply connected
CMC $\frac12$ surfaces in $\h^2\times\R$.
\end{cor}

The fact that $\theta=\frac\pi2$ suggests that this correspondence
looks like the conjugate cousin correspondence between
minimal surfaces in $\R^3$ and CMC $1$ surfaces in $\h^3$
(\cite{bryant}, \cite{umehara}). This correpondence
has nice geometric properties, and is useful to construct
CMC $1$ surfaces in $\h^3$ with some prescribed
geometric properties starting from a solution of a Plateau
problem in $\R^3$ (see for example
\cite{karcher}, \cite{troisdroites}).

In particular, if a minimal surface $\Sigma_1$ in $\nil$ contains an
ambient geodesic $\gamma$, then the normal curvature of $\gamma$
vanishes, and so
$$0=\langle\gamma',\rmS_1\gamma'\rangle
=\langle\gamma',-\rmJ\rmS_2\gamma'+\frac12\rmJ\gamma'\rangle
=-\langle\gamma',\rmJ\rmS_2\gamma'\rangle.$$
This means that $\rmS\gamma'$ is colinear to $\gamma'$, i.e., 
$\gamma$ is a geodesic line of curvature in the
sister CMC $\frac12$ surface in $\h^2\times\R$.

We describe two examples of sister CMC $\frac12$ surfaces in $\h^2\times\R$
of minimal surfaces in $\nil$. We will use the exponential coordinates
given in section \ref{heisenberg} (with $\tau=\frac12$).
We will denote between parentheses ( )
the coordinates of a vector in the coordinate frame 
$(\partial_x,\partial_y,\partial_z)$, and  between brackets [ ]
the coordinates of a vector in the canonical frame $(E_1,E_2,E_3)$;
with these notations one has
$$\left(\begin{array}{c}
a \\
b \\
c
\end{array}\right)=\left[\begin{array}{c}
a \\
b \\
\frac12(ya-xb)+c
\end{array}\right].$$

\begin{example}[vertical plane]
A vertical plane $\cP$ in $\nil$ is a flat minimal surface (but not totally
geodesic). A conformal parametrisation is
$$\varphi:(u,v)\mapsto\left(\begin{array}{c}
v \\
0 \\
u
\end{array}\right).$$
We have $$\varphi_u=E_3,\quad\varphi_v=E_1,\quad
N=E_2,$$
and so $$\nu=0,$$
$$\langle T,\partial_u\rangle=\langle\xi,\varphi_u\rangle
=1,$$
$$\langle T,\partial_v\rangle=\langle\xi,\varphi_v\rangle
=0,$$
i.e., $$T=\partial_u.$$
We also have 
$$\bar\nabla_{\varphi_u}N=\frac12E_1=\frac12\varphi_u,\quad
\bar\nabla_{\varphi_v}N=\frac12E_3=\frac12\varphi_v,$$
so in the direct orthonormal frame $(\partial_u,\partial_v)$ 
we have $$\rmS=-\frac12\left(\begin{array}{cc}
0 & 1 \\
1 & 0
\end{array}\right).$$

We now show that the CMC $\frac12$ sister in $\h^2\times\R$ of $\cP$
is the product $\cH\times\R$ where
$\cH$ is a horocycle in $\h^2$. We will use the upper half-plane model
for $\h^2$. Then $\h^2\times\R=\{(x,y,z)\in\R^3;y>0\}$ and the
metric is $\rmd s^2=\frac1{y^2}(\rmd x^2+\rmd y^2)+\rmd z^2$. We consider
the direct orthonormal frame $(E_1,E_2,E_3)$ defined by
$E_1=y\partial_x$, $E_2=y\partial_y$, $E_3=\partial_z$; it satisfies
$\bar\nabla_{E_1}E_1=E_2$, $\bar\nabla_{E_1}E_2=-E_1$, and the
other derivatives vanish. For $\cH$, we can choose the curve of
equation $y=1$ in $\h^2$. A conformal parametrization of $\cH\times\R$
is
$$\tilde\varphi:(u,v)\mapsto\left(\begin{array}{c}
-u \\
1 \\
v
\end{array}\right).$$
We have $$\tilde\varphi_u=-E_1,\quad\tilde\varphi_v=E_3,\quad
N=E_2,$$
and so $$\tilde\nu=0,\quad\tilde T=\partial_v.$$
We also have 
$$\bar\nabla_{\tilde\varphi_u}N=E_1=-\tilde\varphi_u,\quad
\bar\nabla_{\tilde\varphi_v}N=0,$$
so in the direct orthonormal frame $(\partial_u,\partial_v)$ 
we have $$\tilde\rmS=\left(\begin{array}{cc}
1 & 0 \\
0 & 0
\end{array}\right).$$

Hence, $\tilde\varphi$ induces on $\R^2$ the same metric as $\varphi$,
and we have $\tilde\nu=\nu$, $\tilde T=\rmJ T$ and 
$\tilde\rmS=\rmJ\rmS+\frac12\rmI$, so
$\tilde\varphi$ is the sister immersion of $\varphi$. 
The vertical lines in $\cP$ are mapped to horizontal horocycles
in $\cH\times\R$, and horizontal lines in $\cP$ are mapped to 
vertical lines in $\cH\times\R$.
\end{example}

\begin{example}[surface of equation $z=0$]
The surface $\cA$ of equation $z=0$ in the exponential coordinates is
a minimal surface in $\nil$ which is invariant by rotation about the
$z$-axis (but it is not invariant by any translation; see
\cite{mercuri}). We consider
the following parametrisation:
$$\varphi:(u,v)\mapsto\left(\begin{array}{c}
u\cos v \\
u\sin v \\
0
\end{array}\right),$$
for $u>0$ (the origin in $\cA$ is excluded).
We have $$\varphi_u=\left(\begin{array}{c}
\cos v \\
\sin v \\
0
\end{array}\right)=\left[\begin{array}{c}
\cos v \\
\sin v \\
0
\end{array}\right],$$
$$\varphi_v=\left(\begin{array}{c}
-u\sin v \\
u\cos v \\
0
\end{array}\right)=\left[\begin{array}{c}
-u\sin v \\
u\cos v \\
-\frac12u^2
\end{array}\right],$$
so $$\langle\varphi_u,\varphi_u\rangle=1,$$
$$\langle\varphi_v,\varphi_v\rangle=u^2\left(1+\frac{u^2}4\right),$$
$$\langle\varphi_u,\varphi_v\rangle=0.$$
The unit normal vector is 
$N=\frac{\varphi_u\times\varphi_v}{||\varphi_u\times\varphi_v||}$;
we compute that $$\nu=\frac1{\sqrt{1+\frac{u^2}4}}.$$
A direct orthonormal frame $(e_1,e_2)$ is given by
$$e_1=\partial_u,\quad e_2=\frac1{u\sqrt{1+\frac{u^2}4}}\partial_v.$$
We compute that $$T=-\frac u{2\sqrt{1+\frac{u^2}4}}\partial_v.$$

We now show that the CMC $\frac12$ sister in $\h^2\times\R$ of $\cA$
is the CMC $\frac12$ graph $\cB$ of theorem D in \cite{nellicmc}. This
surface $\cB$ is also invariant by rotation about a vertical axis.
If we take for $\h^2$ the Poincar\'e unit disk model, then $\cB$ is
the graph of the function $(x,y)\mapsto\frac2{\sqrt{1-x^2-y^2}}$.
We will use the Lorentzian for $\h^2\times\R$, i.e., 
$$\h^2\times\R=\{(x^0,x^1,x^2,x^3)\in\LL^3\times\R;
-(x^0)^2+(x^1)^2+(x^2)^2=-1,x_0>0\}$$ with the restriction of the
quadratic form $-(\rmd x^0)^2+(\rmd x^1)^2
+(\rmd x^2)^2+(\rmd x^3)^2$. In this model, we consider
the map
$$\tilde\varphi:(u,v)\mapsto\left(\begin{array}{c}
1+\frac{u^2}2 \\
u\sqrt{1+\frac{u^2}4}\cos v \\
u\sqrt{1+\frac{u^2}4}\sin v \\
2\sqrt{1+\frac{u^2}4}
\end{array}\right),$$
for $u>0$. We can check that it is a parametrization of $\cB$ minus
the origin (using that the correspondence between the Poincar\'e model
and the Lorentzian model is given by $x+iy=\frac{x^1+ix^2}{1+x^0}$,
$z=x^3$). We have 
$$\tilde\varphi_u=\frac1{\sqrt{1+\frac{u^2}4}}
\left(\begin{array}{c}
u\sqrt{1+\frac{u^2}4} \\
1+\frac{u^2}2\cos v \\
1+\frac{u^2}2\sin v \\
\frac u2
\end{array}\right),\quad
\tilde\varphi_v=\left(\begin{array}{c}
0 \\
-u\sqrt{1+\frac{u^2}4}\sin v \\
u\sqrt{1+\frac{u^2}4}\cos v \\
0
\end{array}\right),$$
so $$\langle\tilde\varphi_u,\tilde\varphi_u\rangle=1,$$
$$\langle\tilde\varphi_v,\tilde\varphi_v\rangle
=u^2\left(1+\frac{u^2}4\right),$$
$$\langle\tilde\varphi_u,\tilde\varphi_v\rangle=0,$$
so $\tilde\varphi$ induces the same metric as $\varphi$.
We compute that $$\tilde T=\frac u{2\sqrt{1+\frac{u^2}4}}e_1
=\rmJ T.$$
Thus we also have $\tilde\nu^2=\nu^2$. Moreover, $\tilde\varphi_u$
points outwards and $\tilde\varphi_v$ 
points in the counter-clockwise direction, so the normal $\tilde N$
points up, i.e., $\tilde\nu>0$. So we get $$\tilde\nu=\nu.$$
It remains to check that $\tilde\rmS=\rmJ\rmS+\frac12\rmI$. 
Since $\nu\neq 0$, the 
compatibility equations \eqref{conditionT1} for $\varphi$ and
$\tilde\varphi$ imply that
$\tilde\rmS=\rmJ(\rmS-\frac12\rmJ)=\rmJ\rmS+\frac12\rmI$. 
Hence $\tilde\varphi$ is the sister immersion of $\varphi$.

The straight lines in $\cA$ passing through the origin are mapped 
to the generatrices of $\cB$, which are lines of curvatures lying in
vertical planes. Thus the symmetries of $\cB$ with respect to these
vertical planes correspond to the symmetries of $\cA$ with respect
to the straight lines passing through the origin.
\end{example}

\begin{example}[CMC rotational spheres]
The sister of the CMC $H_1$ rotational sphere in $\nil$ is the
CMC $\sqrt{H_1^2+\frac14}$ rotational sphere in $\h^2\times\R$. Indeed,
the sister of this sphere is a possibly immersed CMC
sphere in $\h^2\times\R$, which is necessarily rotational by a
theorem of Abresch and Rosenberg (\cite{abresch}).
\end{example}

\begin{rem}
CMC $H$ surfaces in $\h^2\times\R$ have very different properties
when $H\leqslant\frac12$ and when $H>\frac12$; for example compact 
embedded CMC
$H$ surfaces exist only for $H>\frac12$. The reader can refer for 
example to \cite{nellicmc}. An explanation is that CMC $H$ surfaces
in $\h^2\times\R$ arise from minimal surfaces in a Berger sphere when
$H>\frac12$, in $\nil$ when $H=\frac12$, and in a space $\psl$ when
$H<\frac12$.
\end{rem}

\begin{rem}
When $\kappa-4\tau^2=0$, the sister relation is the composition
of the classical cousin relation between the round $3$-spheres and $\R^3$
and of the conjugation by a phase $\theta$ in the associate
family. The hyperbolic $3$-space
does not appear in this classification since it is not a fibration
over a $2$-manifold of constant curvature.
\end{rem}

\begin{rem}
A classical problem in the theory of minimal surfaces is the
question of the existence of minimal isometric deformations of 
a given minimal surface.
The compatibility equations show that an associated family of a
given minimal surface (i.e., a one-parameter family of minimal isometric
deformation of this surface obtained by rotating the shape operator)
in a homogeneous $3$-manifold
$\E$ when $\tau\neq 0$ cannot be obtained in a simple way as in
$\s^3$, $\R^3$, $\h^3$, $\s^2\times\R$ or $\h^2\times\R$ 
(see \cite{codazzi}). Indeed, if the quadruple
$(\rmd s^2,\rmS,T,\nu)$ satisfies
the compatibility equations for $\E$, then, in general, the quadruple
$(\rmd s^2,e^{\theta\rmJ}\rmS,e^{\theta\rmJ}T,\nu)$ where 
$\theta\in\R\setminus2\pi\Z$ does not. The question of the existence
of the associate family for minimal surfaces in $\E$ when $\tau\neq 0$
remains open.
\end{rem}

\subsection{Twin immersions}

In this section we will study the special case of sister immersions lying
in the same homogeneous $3$-manifold. They necessarily have
opposite mean curvatures.

\begin{thm} \label{twins}
Let $\E$ be a homogeneous $3$-manifold with a $4$-dimensional
isometry group, of base curvature $\kappa$ and bundle curvature $\tau$.
Let $\xi$ be its vertical vector field.

Let $\Sigma$ be a simply connected
Riemann surface and let $x:\Sigma\to\E$ be
a conformal constant mean curvature $H\neq 0$ immersion.
Let $N$ be the induced normal (compatible with the orientation
of $\Sigma$). Let
$\rmS$ be the symmetric operator on $\Sigma$ induced by the shape
operator of $x(\Sigma)$ associated to the normal $N$. Let $T$
be the vector field on $\Sigma$ such that $\rmd x(T)$ is the
projection of $\xi$ onto $\rmT(x(\Sigma))$. Let
$\nu=\langle N,\xi\rangle$. Let $$\theta=-2\arctan\frac H\tau.$$

Then there exists a unique conformal immersion
$\hat x:\Sigma\to\E$ such that:
\begin{enumerate}
\item the metrics induced on $\Sigma$ by $x$ and $\hat x$ are the same,
\item the symmetric operator on $\Sigma$ induced by the shape operator
of $\hat x(\Sigma)$ is $\tilde\rmS=e^{\theta\rmJ}(\rmS-H\rmI)-H\rmI
=e^{\theta\rmJ}(\rmS-\tau\rmJ)+\tau\rmJ$,
\item $\xi=\rmd\hat x(e^{\theta\rmJ}T)+\nu\hat N$ where 
$\hat N$ is the unit normal to $\hat x$.
\end{enumerate}

Moreover, this immersion $\hat x$ is unique up to isometries of $\E$
preserving the orientations of both the fibers and the base of the
fibration, and it has constant mean curvature $-H$.

It is called the twin immersion of the immersion $x$.
\end{thm}

\begin{proof}
This is a particular case of theorem \ref{sisters} with
$\E_1=\E_2=\E$, $\tau_1=-\tau_2=\tau$, $H_1=-H_2=H$.
It sufficies to check that the phase $\theta$
satisfies $\tau-iH=e^{i\theta}(\tau+iH)$.

The equivalence of the two expressions of $\tilde\rmS$ is a 
consequence of \eqref{rotationshape2}.
\end{proof}

%

We notice that when $\tau\to 0$, then
$\theta\to\pi$, i.e., $\tilde T\to-T$, and also $\tilde\rmS\to-\rmS$. 
This limit corresponds to
the image of the initial surface by a horizontal symmetry in 
$\M^2(\kappa)\times\R$.

Moreover, we notice that the twin surface of a multigraph (over a part of
the base of the fibration) is also a multigraph
(since a surface is a multigraph if and only if $\nu$ does not vanish).

This suggests that the twin surface could be used to get an
Alexandrov reflection-type principle 
in homogeneous manifolds with non-vanishing bundle curvature, since there is no
Alexandrov reflection principle (see \cite{alexandrov})
in these manifolds (the horizontal and 
vertical ``symmetries'' are not isometries). Such an Alexandrov reflection
principle would be very useful for the theory of CMC surfaces
in homogeneous manifolds, in particular for proving that any closed
embedded CMC surface in the Heisenberg space or in $\psl$
is a rotational sphere (this was proved 
for CMC surfaces in $\R^3$, $\h^3$, a $3$-hemisphere, $\h^2\times\R$ and 
a $2$-hemisphere cross $\R$ using the Alexandrov reflection principle).

We now give some examples of twin surfaces in the Heisenberg space
$\nil$ with its standard metric (i.e., $\kappa=0$, $\tau=\frac12$). We
will use the exponential coordinates described in section 
\ref{heisenberg}. Figueroa, Mercuri and Pedrosa classified CMC surfaces
in $\nil$ invariant by a one-parameter family of translations or 
rotations (see \cite{mercuri}; note that in their article the mean
curvature is defined as the trace of the shape operator, whereas in this 
paper it is defined as the half of the trace). We will compute the
twin surfaces of these examples. 
We will denote between parentheses ( )
the coordinates of a vector in the coordinate frame 
$(\partial_x,\partial_y,\partial_z)$, and  between brackets [ ]
the coordinates of a vector in the canonical frame $(E_1,E_2,E_3)$.

\begin{example}[translational tubes] \label{tube}
Let $H>0$. The map
$$\varphi:(u,v)\mapsto\left(\begin{array}{c}
u \\
\frac{\cos v}{2H} \\
u\frac{\cos v}{4H}+\frac1{4H}f(v)
\end{array}\right),$$ with
$$f(v)=\sqrt{1+\frac{\cos^2v}{4H^2}}\sin v
+\frac{1+4H^2}{2H}\arcsin\left(\frac{\sin v}{\sqrt{1+4H^2}}\right),$$
for $(u,v)\in\R^2$, is a CMC $H$ immersion defining
a surface which is invariant by
horizontal translations in the $x$-direction. This surface is an annulus,
and it is a bigraph over a part of the minimal surface of equation 
$z=\frac{xy}2$; moreover it is ``symmetric'' with respect to this
minimal surface.

We have $$\varphi_u=\left(\begin{array}{c}
1 \\
0 \\
\frac{\cos v}{4H}
\end{array}\right)=\left[\begin{array}{c}
1 \\
0 \\
\frac{\cos v}{2H}
\end{array}\right],$$
$$\varphi_v=\left(\begin{array}{c}
0 \\
-\frac{\sin v}{2H} \\
-u\frac{\sin v}{4H}+\frac1{4H}f'(v)
\end{array}\right)=\left[\begin{array}{c}
0 \\
-\frac{\sin v}{2H} \\
\frac1{4H}f'(v)
\end{array}\right],$$
$$f'(v)=2\cos v\sqrt{1+\frac{\cos^2v}{4H^2}},$$
and so
$$\langle\varphi_u,\varphi_u\rangle=1+\frac{\cos^2v}{4H^2},$$
$$\langle\varphi_v,\varphi_v\rangle
=\frac1{4H^2}\left(1+\frac{\cos^4v}{4H^2}\right).$$
$$\langle\varphi_u,\varphi_v\rangle
=\frac{\cos^2v}{4H^2}\sqrt{1+\frac{\cos^2v}{4H^2}}.$$
The unit normal vector is given by 
$N=\frac{\varphi_u\times\varphi_v}{||\varphi_u\times\varphi_v||}$;
we compute that $$\nu=-\frac{\sin v}
{\sqrt{1+\frac{\cos^4v}{4H^2}}}.$$

We have $$\langle T,\partial_u\rangle
=\langle\xi,\varphi_u\rangle=\frac{\cos v}{2H},$$
$$\langle T,\partial_v\rangle
=\langle\xi,\varphi_v\rangle
=\frac{\cos v}{2H}\sqrt{1+\frac{\cos^2v}{4H^2}},$$

We notice that $\nu(u_1,-v)=-\nu(u_2,v)$ for all $(u_1,u_2,v)$.
This indicates that the twin immersion could be an 
orientation-reversing reparametrization of the surface.
For this reason we set
$$\tilde\varphi:(u,v)\mapsto\varphi(u+h(v),-v)
=\left(\begin{array}{c}
u+h(v) \\
\frac{\cos v}{2H} \\
(u+h(v))\frac{\cos v}{4H}-\frac1{4H}f(v)
\end{array}\right)$$ where $h$ is a function.
This is a CMC $-H$ immersion defining globally
the same surface as $\varphi$.
We compute that $$\tilde\varphi_u=\left[\begin{array}{c}
1 \\
0 \\
\frac{\cos v}{2H}
\end{array}\right],\quad
\tilde\varphi_v=\left[\begin{array}{c}
h'(v) \\
-\frac{\sin v}{2H} \\
h'(v)\frac{\cos v}{2H}-\frac1{4H}f'(v)
\end{array}\right],$$
and so
$$\langle\tilde\varphi_u,\tilde\varphi_u\rangle
=1+\frac{\cos^2v}{4H^2},$$
\begin{eqnarray*}
\langle\tilde\varphi_v,\tilde\varphi_v\rangle & = &
\left(1+\frac{\cos^2v}{4H^2}\right)h'(v)^2
-\frac{\cos^2v}{2H^2}h'(v)\sqrt{1+\frac{\cos^2v}{4H^2}} \\
& & +\frac1{4H^2}\left(1+\frac{\cos^4v}{4H^2}\right),
\end{eqnarray*}
$$\langle\tilde\varphi_u,\tilde\varphi_v\rangle
=\left(1+\frac{\cos^2v}{4H^2}\right)h'(v)
-\frac{\cos^2v}{4H^2}\sqrt{1+\frac{\cos^2v}{4H^2}}.$$
Thus $\tilde\varphi$ induces on $\R^2$ the same metric as $\varphi$ if
and only if 
$$h'(v)=\frac{\cos^2v}{2H^2\sqrt{1+\frac{\cos^2v}{4H^2}}}.$$

We now assume that this condition is satisfied; we can also assume that
$h(0)=0$. The function $h$ is increasing. We have 
$$\tilde\nu=\nu,$$
$$\langle\tilde T,\partial_u\rangle
=\langle\xi,\tilde\varphi_u\rangle=\frac{\cos v}{2H},$$
$$\langle\tilde T,\partial_v\rangle
=\langle\xi,\tilde\varphi_v\rangle
=\frac{\cos v}{2H\sqrt{1+\frac{\cos^2v}{4H^2}}}
\left(\frac{\cos^2v}{4H^2}-1\right).$$

The direct orthonormal frame $(e_1,e_2)$ obtained from the frame 
$(\partial_u,\partial_v)$ by the Gram-Schmidt process satisfies
$$e_1=\frac{\partial_u}{||\partial_u||},$$
$$e_2=
\frac{-\langle\partial_u,\partial_v\rangle\partial_u
+||\partial_u||^2\partial_v}
{||\partial_u||
\sqrt{||\partial_u||^2||\partial_u||^2
-\langle\partial_u,\partial_v\rangle^2}}.$$
A computation gives
$$||\partial_u||^2||\partial_u||^2
-\langle\partial_u,\partial_v\rangle^2
=\frac1{4H^2}\left(1+\frac{\cos^2v}{4H^2}\right).$$
Thus we get
$$e_1=\frac1{\sqrt{1+\frac{\cos^2v}{4H^2}}}\partial_u,$$
$$e_2=-\frac{\cos^2v}{2H\sqrt{1+\frac{\cos^2v}{4H^2}}}\partial_u
+2H\partial_v.$$

So we have $$T=\frac{\cos v}{\sqrt{1+\frac{\cos^2v}{4H^2}}}
\left(\frac1{2H}e_1+e_2\right),$$
$$\tilde T=\frac{\cos v}{\sqrt{1+\frac{\cos^2v}{4H^2}}}
\left(\frac1{2H}e_1-e_2\right).$$

Let $\theta=-2\arctan(2H)$. Then we have
$$\cos\theta=\frac{1-4H^2}{1+4H^2},\quad
\cos\theta=-\frac{4H}{1+4H^2}.$$
Since $\rmJ e_1=e_2$ and $\rmJ e_2=-e_1$, we get
$$e^{\theta\rmJ}T=\tilde T.$$
Finally, the compatibility equation \eqref{conditionT1} implies that
$$\tilde S=e^{\theta\rmJ}(\rmS-\tau\rmJ)+\tau\rmJ$$
at points where $\nu\neq 0$; and by continuity this identity
holds everywhere. This proves that $\tilde\varphi$ is the twin
immersion of $\varphi$.

Thus the translational tube is \emph{globally} invariant by the
twin relation, but it is \emph{not pointwise} invariant: the
correspondence is
$$\varphi(u,v)\mapsto\varphi(u+h(v),-v).$$
Geometrically, this correspondence maps a point of the tube to the
other point of the tube lying in the same fiber
and then translates it by $h(v)$ in the $x$-direction.
In particular, the closed curve $v\mapsto\varphi(u_0,v)$ is
mapped to the curve $v\mapsto\varphi(u_0+h(v),-v)$,
which is \emph{not} closed.
\end{example}

\begin{example}[rotational spheres] \label{sphere}
Let $H>0$. The map
$$\varphi:(u,v)\mapsto\left(\begin{array}{c}
\frac1H\cos u\cos v \\
\frac1H\sin u\cos v \\
\frac1{2H}f(v)
\end{array}\right),$$ with
$f$ as in example \ref{tube},
for $(u,v)\in\R\times(-\frac{\pi}2,\frac{\pi}2),$
is a CMC $-H$ immersion defining a rotational sphere minus the top and
bottom points (the normal of the immersion points outside whereas the
mean curvature vector points inside).
It is a bigraph over a part of the minimal surface of equation 
$z=0$; moreover it is ``symmetric'' with respect to this
minimal surface.

We have $$\varphi_u=\frac1H\left[\begin{array}{c}
-\sin u\cos v \\
\cos u\cos v \\
-\frac1{2H}\cos^2v
\end{array}\right],\quad
\varphi_v=\frac1H\left[\begin{array}{c}
-\cos u\sin v \\
-\sin u\sin v \\
\frac1{2}f'(v)
\end{array}\right],$$
and so
$$\langle\varphi_u,\varphi_u\rangle
=\frac{\cos^2v}{H^2}\left(1+\frac{\cos^2v}{4H^2}\right),$$
$$\langle\varphi_v,\varphi_v\rangle
=\frac1{H^2}\left(1+\frac{\cos^2v}{4H^2}\right),$$
$$\langle\varphi_u,\varphi_v\rangle
=-\frac{\cos^3v}{2H^3}\sqrt{1+\frac{\cos^2v}{4H^2}}.$$
The unit normal vector is given by 
$N=\frac{\varphi_u\times\varphi_v}{||\varphi_u\times\varphi_v||}$;
we compute that $$\nu=\frac{\sin v}
{\sqrt{1+\frac{\cos^4v}{4H^2}}}.$$

We have $$\langle T,\partial_u\rangle
=\langle\xi,\varphi_u\rangle=-\frac{\cos^2v}{2H^2},$$
$$\langle T,\partial_v\rangle
=\langle\xi,\varphi_v\rangle
=\frac{\cos v}{H}\sqrt{1+\frac{\cos^2v}{4H^2}}.$$

Let $$\tilde\varphi:(u,v)\mapsto\varphi(u+g(v),-v)
=\left(\begin{array}{c}
\frac1H\cos(u+g(v))\cos v \\
\frac1H\sin(u+g(v))\cos v \\
-\frac1{2H}f(v)
\end{array}\right)$$
where $g$ is a function. This is a CMC $H$ immersion defining globally
the same surface as $\varphi$.
We compute that $$\tilde\varphi_u=\frac1H\left[\begin{array}{c}
-\sin(u+g(v))\cos v \\
\cos(u+g(v))\cos v \\
-\frac1{2H}\cos^2v
\end{array}\right],$$
$$\tilde\varphi_v=\frac1H\left[\begin{array}{c}
-\cos(u+g(v))\sin v-g'(v)\sin(u+g(v))\cos v \\
-\sin(u+g(v))\sin v+g'(v)\cos(u+g(v))\cos v \\
-\frac12f'(v)-\frac1{2H}g'(v)\cos^2v
\end{array}\right],$$
and thus
$\tilde\varphi$ induces on $\R\times(\frac{\pi}2,\frac{\pi}2)$
the same metric as $\varphi$ if
and only if 
$$g'(v)=-\frac{\cos v}{H\sqrt{1+\frac{\cos^2v}{4H^2}}}.$$

We now assume that this condition is satisfied; we can also assume that
$g(0)=0$. The function $g$ is odd and $2\pi$-periodic. We have 
$$\tilde\nu=\nu,$$
$$\langle\tilde T,\partial_u\rangle
=\langle\xi,\tilde\varphi_u\rangle=-\frac{\cos^2v}{2H^2},$$
$$\langle\tilde T,\partial_v\rangle
=\langle\xi,\tilde\varphi_v\rangle
=\frac{\cos v}{H\sqrt{1+\frac{\cos^2v}{4H^2}}}
\left(\frac{\cos^2v}{4H^2}-1\right).$$

The direct orthonormal frame $(e_1,e_2)$ obtained from the frame 
$(\partial_u,\partial_v)$ by the Gram-Schmidt process satisfies
$$e_1=\frac H{\cos v\sqrt{1+\frac{\cos^2v}{4H^2}}}\partial_u,$$
$$e_2=-\frac{\cos v}{2\sqrt{1+\frac{\cos^2v}{4H^2}}}\partial_u
+H\partial_v.$$

So we have $$T=\frac{\cos v}{\sqrt{1+\frac{\cos^2v}{4H^2}}}
\left(-\frac1{2H}e_1+e_2\right),$$
$$\tilde T=\frac{\cos v}{\sqrt{1+\frac{\cos^2v}{4H^2}}}
\left(-\frac1{2H}e_1-e_2\right).$$

Let $\theta=2\arctan(2H)$. We check as in example \ref{tube}
that $$e^{\theta\rmJ}T=\tilde T,$$
$$\tilde S=e^{\theta\rmJ}(\rmS-\tau\rmJ)+\tau\rmJ.$$
This proves that $\tilde\varphi$ is the twin
immersion of $\varphi$.

Thus the rotational sphere is \emph{globally} invariant by the
twin relation, but it is \emph{not pointwise} invariant: the
correspondence is
$$\varphi(u,v)\mapsto\varphi(u+g(v),-v).$$
Geometrically, this correspondence maps a point of the sphere to the
other point of the sphere lying in the same fiber
and then rotates it by the angle $g(v)$ about the $z$-axis.
In particular, the circle $v\mapsto\varphi(u_0,v)$ lying in a
vertical plane is
mapped to the curve $v\mapsto\varphi(u_0+g(v),-v)$,
which is closed but not contained in a vertical plane.
\end{example}

\bibliographystyle{alpha}
\bibliography{immersions}

\end{document}